\documentclass{svmult}
\usepackage{amsmath,amsfonts,amssymb}

\spnewtheorem{theorem*}{Theorem}[section]{\bf}{\it} 
\spnewtheorem{lemma*}[theorem*]{Lemma}{\bf}{\it}
\spnewtheorem{proposition*}[theorem*]{Proposition}{\bf}{\it}
\spnewtheorem{corollary*}[theorem*]{Corollary}{\bf}{\it}
\spnewtheorem{definition*}[theorem*]{Definition}{\bf}{\it}
\spnewtheorem{example*}[theorem*]{Example}{\bf}{\it} 
\spnewtheorem{remark*}[theorem*]{Remark}{\bf}{\it} 

\def\a{\mathfrak{a}}
\def\al{\alpha}
\def\A{{\rm A}}
\def\AMN{\A_N^M}

\def\B{{\rm B}}

\def\bh{{\bar h}}
\def\bi{{\bar\imath}}
\def\bj{{\bar\jmath}}
\def\bk{{\bar k}}
\def\bl{{\bar l}}
\def\Bp{{\rm B}'}
\def\Bpp{{\rm B}''}

\def\C{{\rm C}}
\def\CC{{\mathbb C}}
\def\CM{\CC^{M|M}}
\def\CN{\CC^{N|N}}
\def\CNM{\CC^{\ts N+M\ts|\ts N+M}}
\def\com{\ts,\hskip-.5pt}

\def\D{{\rm D}}
\def\de{\delta}

\def\deg{{{\rm deg}\ts}}

\def\End{{\rm End}\hskip1pt}
\def\EndCN{\End(\CN)}

\def\et{\eta\ts}

\def\g{\mathfrak{g}}
\def\ge{\geqslant}
\def\glM{\mathfrak{gl}_{M|M}}
\def\glN{\mathfrak{gl}_{N|N}}
\def\glNM{\mathfrak{gl}_{\ts N+M\ts|\ts N+M}}

\def\I{{\rm I}}
\def\id{{{\rm id}}}

\def\J{{\rm J}}

\def\lc{{\ts,\hskip.95pt\ldots\ts,\ts}}
\def\le{\leqslant}

\def\O{\mathcal{O}}
\def\om{\omega}
\def\op{\oplus}
\def\ot{\otimes}

\def\ph{\varphi}
\def\ps{\psi}

\def\q{\mathfrak{q}}
\def\qM{\q_M}
\def\qN{\q_N}
\def\qNM{\q_{\ts N+M}}
\def\Q{{\rm Q}}

\def\S{{\rm S}}
\def\si{\sigma}

\def\ts{{\hskip1pt}}

\def\U{{\rm U}}
\def\Ug{\U(\g)}
\def\UM{{\U}(\qM)}
\def\UNM{{\U}(\qNM)}
\def\UN{{\U}(\qN)}

\def\Vb{\,\overline{\hskip-1pt V\hskip-1pt}\,}

\def\X{{\rm X}}
\def\Xt{\widetilde{X}}

\def\YN{{\rm Y}(\qN)}
\def\Yt{\widetilde{Y}}

\def\ZN{{\rm Z}(\qN)}
\def\ZNM{{\rm Z}(\qNM)}

\def\ZZ{\mathbb Z}

\begin{document} 
 
\titlerunning{Yangian of the queer Lie superalgebra}
\title*{Centralizer construction of the Yangian
of the queer Lie superalgebra}
\author{Maxim Nazarov\inst{1}
\and Alexander Sergeev\inst{2}}
\institute{
Department of Mathematics, 
University of York, 
Heslington,
York YO10 5DD, 
England; 
\texttt{mln1@york.ac.uk}
\and 
Department of Mathematics, 
Balakovo Institute of Technology, 
Balakovo 413800,
Russia;
\texttt{sergeev@bittu.org.ru}
}

\maketitle

                                                             
\renewcommand{\theequation}{\thesection.\arabic{equation}} 
\makeatletter                                          
\@addtoreset{equation}{section}                        
\makeatother                                           
                                                       

\medskip\noindent 
\textit{To Professor Anthony Joseph on the occasion of his 60\ts th birthday} 


\begin{abstract} 
Consider the complex matrix Lie superalgebra $\glN$ with the
standard generators $E_{ij}$ where $i\com j=\pm1\lc\!\pm N$.
Define an involutive automorphism $\et$ of $\glN$ by
$\et(E_{ij})=E_{-i,-j}$. The queer Lie superalgebra $\qN$ is the fixed
point subalgebra in $\glN$ relative to $\et$.
Consider the twisted polynomial current Lie superalgebra
$$
\g=\{\,X(t)\in\glN[t]\,:\,\et(X(t))=X(-t)\,\}\,.
$$
The enveloping algebra $\Ug$ of the Lie superalgebra 
$\g$ has a deformation,
called the Yangian of $\qN$. For each $M=1\com2\com\,\ldots$
denote by $\AMN$ the centralizer
of $\qM\subset\qNM$
in the associative superalgebra $\UNM\ts$.
We construct a sequence of surjective homomorphisms
$
\UN\leftarrow
{\rm A}_N^1
\leftarrow
{\rm A}_N^2
\leftarrow\ldots\,.
$
We describe the inverse limit of the sequence of centralizer algebras
${\rm A}_N^1,{\rm A}_N^2,\,\ldots$ in terms of the Yangian of $\qN$.
\end{abstract}


\section{\hskip-5pt.\hskip6pt Main results} 

In this article we work with the queer Lie superalgebra $\qN$.
This is perhaps the most interesting super-analogue of the general linear Lie
algebra $\mathfrak{gl}_N$, see for instance \cite{S2}. We will realize $\qN$
as a subalgebra in the general linear Lie superalgebra $\glN$
over the complex field $\CC$.
Let the indices $i\com j$ run through $-N\lc\!\!-\!1\com1\lc N$.
Put $\bi=0$ if $i>0$ and $\bi=1$ if $i<0$. Take the
$\ZZ_2$\ts-\ts graded vector space $\CN$.
Let $e_i\in\CN$ be the standard basis vectors.
The $\ZZ_2$\ts-\ts gradation on $\CN$ is defined so that $\deg e_i=\bi$.
Let $E_{ij}\in\EndCN$ be the matrix units\ts: $E_{ij}\ts e_k=\de_{jk}\ts e_i$.
The algebra $\EndCN$ is $\ZZ_2$\ts-\ts graded so that 
$\deg E_{ij}=\bi\ts+\bj$.
We will also regard $E_{ij}$ as basis elements of the Lie superalgebra $\glN$.
The {\it queer\/} Lie superalgebra $\qN$ is
the fixed point subalgebra in $\glN$ with respect to
the involutive automorphism $\eta$ defined by
\begin{equation}\label{1.1}
\et:\thinspace E_{ij}\mapsto E_{-i,-j}\,.
\end{equation}
Thus as a vector subspace, $\qN\subset\glN$ is spanned by the elements
$$
F_{ij}=E_{ij}+E_{-i,-j}\,.
$$ 
Note that $F_{-i,-j}=F_{ij}\,$.
The elements $F_{ij}$ with $i>0$ form a basis of $\qN$.

The vector subspace of $\End(\CN)$ spanned by the elements $F_{ij}$
is closed with respect to the usual matrix multiplication. Hence we can 
also regard it as an associative algebra. Denote this associative
algebra by $\Q_N$, to distinguish its structure from that of
the Lie superalgebra $\qN$. Both $\EndCN$ and $\Q_N$
are simple as associative $\ZZ_2$\ts-\ts graded algebras, see
\cite[Theorem 2.6]{J}.

The enveloping algebra $\UN$ of the Lie superalgebra $\qN$ is
a $\ZZ_2$\ts-\ts graded associative unital algebra. In this article
we will always keep to the following convention.
Let $\A$ and $\B$ be any two associative
$\ZZ_2$\ts-\ts graded algebras. Their tensor product
$\A\ot\B$ is
a $\ZZ_2$\ts-\ts graded algebra such that for any homogeneous
elements $X,X'\in\A$ and $Y,Y'\in\B$
\begin{eqnarray}\label{otp}
(X\ot Y)\ts (X'\ot Y')
&=&
(-1)^{\,\deg X'\deg Y}\,
X\ts X'\ot Y\ts Y',
\\[4pt]
\label{otd}
\deg\ts(X\ot Y)&=&\deg X+\deg Y\ts.
\end{eqnarray}

By definition, an anti-homomorphism $\om:\A\to\B$ is any linear map
which preserves the $\ZZ_2\ts$-\ts gradation and satisfies
any homogeneous $X,X'\in\A$
\begin{equation}\label{ai}
\om(\ts X\ts X')=(-1)^{\,\deg X\,\deg X'}\,\om(X')\,\,\om(X)\,.
\end{equation}
For any Lie superalgebra $\a$, the \textit{principal\/}
anti-automorphism of the enveloping 
$\ZZ_2\ts$-\ts graded algebra $\U(\a)$ is
determined by the assignment $X\mapsto-X$ for $X\in\a\ts$.

The {\it supercommutator\/} of any two homogeneous elements
$X,Y\in\A$ is by definition
\begin{equation}\label{scom}
[\ts X\com Y\ts]\,=\,XY-
(-1)^{\,\deg X\deg Y}\,YX.
\vspace{2pt}
\end{equation}
This definition extends to arbitrary elements $X,Y\in\A$ by linearity. 
It is the bracket (\ref{scom}) that defines the Lie superalgebra
structure on the vector space $\A=\EndCN\ts$.
Thus for any indices $i\com j\com k\com l=\pm1\lc\!\pm\!N$ we have
$$
\begin{aligned}[t]
[\ts E_{ij}\com E_{kl}\ts]
&\,=\,
\de_{kj}\,E_{il}-
(-1)^{\,(\ts\bi\,+\ts\bj\ts)\ts(\ts\bk\,+\,\bl\ts)}\,\de_{il}\,E_{kj}\,;
\\[3pt]
[\ts F_{ij}\com F_{kl}\ts]
&\,=\,
\de_{kj}\,F_{il}-
(-1)^{\,(\ts\bi\,+\ts\bj\ts)\ts(\ts\bk\,+\,\bl\ts)}\,\de_{il}\,F_{kj}\,\,+
\end{aligned}
$$
\vskip-7pt
\begin{equation}
\label{FR}
\de_{-k,j}\,F_{-i,l}-
(-1)^{\,(\ts\bi\,+\ts\bj\ts)\ts(\ts\bk\,+\,\bl\ts)}\,
\de_{i,-l}\,F_{k,-j}\,.
\end{equation}

For any $\A$ and any subset $\C\subset\A$,
by the {\it centralizer\/} of $\C$ in ${\rm A}$
we mean the collection of all elements $X\in\A$
such that $[\ts X\com Y\ts]=0$ for any $Y\in\C$.
To remind the reader about this convention, we shall then
refer to the $\ZZ_2$\ts-\ts graded algebra $\A$ as a superalgebra.
The centre $\ZN$ of the enveloping algebra $\UN$ will be always taken
in the superalgebra sense. A set of generators of the algebra
$\ZN$ was given in \cite{S1}. In particular, all central elements
of $\UN$ were shown to have $\ZZ_2$\ts-\ts degree 0.
A distinguished basis of the vector space $\ZN$ was constructed in \cite{N2}.

Let us recall the principal results of \cite{S1} here.
For any indices $n\ge1$ and $i\com j=\pm1\lc\!\pm\!N$
denote by $F_{ij\ts|\ts N}^{\,(n)}$ the element of the algebra $\UN$
\begin{equation}\label{FN}
\sum_{k_1,\ldots,\ts k_{n-1}}
(-1)^{\,\bk_1\,+\,\ldots\,+\,\bk_{n-1}}\,
F_{i\ts k_1} F_{k_1k_2}
\ts\ldots\, 
F_{k_{n-2}k_{n-1}} F_{k_{n-1}j}
\end{equation}
where each of the indices $k_1\lc k_{n-1}$
runs through $\pm1\lc\pm N$. Note that 
\begin{equation}\label{FN-}
F_{-i,-j\ts|\ts N}^{\,(n)}=
(-1)^{\,n-1}\,
F_{ij\ts|\ts N}^{\,(n)}\,.
\end{equation}
Of course, here $F_{ij\ts|\ts N}^{\,(1)}=F_{ij}\ts$.
Observe that if $n>1$, then by the definition~(\ref{FN})
\begin{equation}\label{split}
F_{ij\ts|\ts N}^{\,(n)}\,=\,\sum_k\,
(-1)^{\ts\bk}\,
F_{ik}\,F_{kj\ts|\ts N}^{\,(n-1)}
\end{equation}
where the index $k$ runs through $\pm1\lc\pm N$.
Using this observation one proves by induction on $n=1\com2\com\ts\ldots$
the following generalization of (\ref{FR})\ts:
in the $\ZZ_2\ts$-\ts graded algebra $\UN$ the supercommutator
$$
[\ts F_{ij}\com F_{kl\ts|\ts N}^{\,(n)}\ts]\,=\,
\de_{kj}\,F_{il\ts|\ts N}^{\,(n)}-
(-1)^{\,(\ts\bi\,+\ts\bj\ts)\ts(\ts\bk\,+\,\bl\ts)}\,
\de_{il}\,F_{kj\ts|\ts N}^{\,(n)}\,\,+
$$
\vskip-12pt
\begin{equation}
\label{FNR}
\de_{-k,j}\,F_{-i,l\ts|\ts N}^{\,(n)}-
(-1)^{\,(\ts\bi\,+\ts\bj\ts)\ts(\ts\bk\,+\,\bl\ts)}\,
\de_{i,-l}\,F_{k,-j\ts|\ts N}^{\,(n)}\,\ts.
\end{equation}
For a more general formula, expressing
the supercommutator 
$[\ts F_{ij\ts|\ts N}^{\,(m)}\com F_{kl\ts|\ts N}^{\,(n)}\ts]$
for any $m$ and $n$, see
Proposition \ref{P3.1} and
the remark after its proof.
Now put
\begin{equation}
\label{CN}
C_N^{\,(n)}=
\sum_k\,F_{kk\ts|\ts N}^{\ts(n)}
\end{equation}
where the index $k$ runs through $\pm1\lc\pm N$. The relations
(\ref{FNR}) immediately imply that $C_N^{\,(n)}\in\ZN\ts$.
Note that $C_N^{\,(2)}=C_N^{\,(4)}=\ldots=0$ due to (\ref{FN-}). 
The following proposition has been stated in \cite{S1}
without proof.

\begin{proposition*}\label{P1.0}
The elements $C_N^{\,(1)},C_N^{\,(3)},\,\ldots$
generate the centre\/ $\ZN\ts$.
\end{proposition*}

The dependence of the elements $C_N^{\,(n)}$ and
$F_{ij\ts|\ts N}^{\,(n)}$ of $\UN$ on the index $N$ has been indicated
for the purposes of the next argument, which extends \cite{S1}.

For any integers $N\ge0$ and $M\ge1$ consider the Lie superalgebra $\qNM\ts$.
Now let the indices $i\com j$ run through $-N-M\lc-1\com1\lc N+M$.
Regard the Lie superalgebras $\qN$ and $\qM$ as the subalgebras of
$\qNM$ spanned by the elements $F_{ij}=F_{ij\ts|\ts N+M}^{\,(1)}$ where
$|i|\com|j|\le N$ and $|i|\com|j|>N$ respectively. Denote by
$\AMN$ the centralizer of $\qM$ in the associative superalgebra $\UNM$.

By definition, the centralizer $\AMN$ contains the centre $\ZNM$ of the
$\UNM\ts$. It also contains the subalgebra $\UN\subset\UNM\ts$. 
Moreover, the relations (\ref{FNR}) imply that
the centralizer $\AMN$ contains the elements
\begin{equation}\label{elements}
F_{ij\ts|\ts N+M}^{\,(1)}\,,\,
F_{ij\ts|\ts N+M}^{\,(2)}\,,\,
\,\ldots\ \ \textrm{where}\ \  
|i|\com|j|\le N.
\end{equation}

\begin{theorem*}\label{T1.1}
The elements
$C_{N+M}^{\,(1)}\,,C_{N+M}^{\,(3)}\,,\,\ldots$
and {\rm(\ref{elements})}
generate $\AMN$.
\end{theorem*}

We prove this theorem in Section~2 of the present article.
In the particular case $N=0$, we will then obtain Proposition \ref{P1.0}.

Now take the Lie superalgebra $\q_{\ts N+M-1}\ts$.
As a subalgebra of $\qNM\ts$, it is spanned by the
elements $F_{ij}$ where
$|i|\com|j|<N+M$. In particular, the
subalgebras $\qN$ and $\q_{\ts M-1}$ of
$\q_{\ts N+M-1}$ are spanned by the elements $F_{ij}$ where
$|i|\com|j|\le N$ and $N<|i|\com|j|<N+M$ respectively. 
The enveloping algebra $\U(\q_{\ts N+M-1})$ and its subalgebra
$\A_N^{M-1}$ will be also regarded as subalgebras 
in the associative algebra $\UNM\ts$.
We assume that $M\ge1$ and $\A_N^{\,0}=\UN\ts$.

Denote by $\I_{\ts N+M}$ the right ideal in the algebra $\UNM$
generated by the elements
\begin{equation}\label{JM}
F_{N+M,\ts\pm1}
\,\lc\ts
F_{N+M,\ts\pm(N+M)}
\,.
\end{equation}

\begin{lemma*}\label{L1.3}
\,a) the intersection\/ $\I_{\ts N+M}\cap\AMN$ 
is a two-sided ideal of\/ $\AMN\,;$
\\[2pt]
b)
there is a decomposition\/ 
$\AMN=\A_N^{M-1}\op(\,\I_{\ts N+M}\cap\AMN\ts)\,.$
\end{lemma*}

We prove this lemma in Section 3.
Using Part (b) of the lemma, denote by $\al_M$ the
projection of $\AMN$ to its direct summand $\A_N^{M-1}\ts$.
Due to Part (a), the map $\al_M:\AMN\to\A_N^{M-1}$
is a homomorphism of associative algebras.
The proof of the next proposition will also be given in Section 3.

\begin{proposition*}\label{P1.4}
For any $n\ge1$ and any $i\com j$ such that $|i|\com|j|\le N$ we have
$$
\al_M(F_{ij\ts|\ts N+M}^{\,(n)})=F_{ij\ts|\ts N+M-1}^{\,(n)}
\quad\text{and}\/\quad
\al_M(C_{N+M}^{\,(n)})=C_{N+M-1}^{\,(n)}\,.
$$
\end{proposition*}

The standard filtration (\ref{filtr})
on the enveloping algebra $\UNM$ defines a filtration
on its subalgebra $\AMN$. By definition, the map
$\al_M:\AMN\to\A_N^{M-1}$ preserves that filtration.
It also preserves the $\ZZ_2\ts$-\ts gradation, inherited from
$\UNM\ts$.
Using the homomorphisms $\al_1\com\al_2\com\,\ldots$
define an algebra $\A_N$ as the inverse limit of the sequence
$\A_N^{\,0}\com\A_N^1\com\A_N^2\com\,\ldots\,$
in the category of associative filtered algebras.
The main result of this article is an
explicit description of the algebra $\A_N$
in terms of generators and relations.

By definition, an element of $\A_N$ is any sequence
of elements $Z_0\com Z_1\com Z_2\com\,\ldots$ of the algebras
$\A_N^{\,0}\com\A_N^1\com\A_N^2\com\,\ldots\,$ respectively,
such that $\al_M(Z_M)=Z_{M-1}$ for each $M\ge1\ts$, and
the filtration degrees of the elements in the sequence 
are bounded from above.
Utilising Proposition 1.4, for any $n=1\com2\com\,\ldots$ 
and any $i\com j=\pm1\lc\!\pm\!N$ define an element 
$F_{ij}^{\ts(n)}\in\A_N$ as the sequence
\begin{equation}\label{fseq}
F_{ij\ts|\ts N}^{\,(n)}
\ts\com\,
F_{ij\ts|\ts N+1}^{\,(n)}
\ts\com\,
F_{ij\ts|\ts N+2}^{\,(n)}
\,\com\,
\ldots\,.
\end{equation}
Further, for any $n=1\com3\com\,\ldots$ 
define an element $C^{\ts(n)}\in\A_N$ as the 
sequence
\begin{equation}\label{cseq}
C_N^{\ts(n)}
\com\,
C_{N+1}^{\ts(n)}
\com\,
C_{N+2}^{\ts(n)}
\,\com\,
\ldots\,.
\end{equation}
The filtration degree of
every element in (\ref{fseq}) and (\ref{cseq}) does not
exceed $n\ts$.

Note that the algebra $\A_N$ is unital, and
comes with a $\ZZ_2\ts$-\ts gradation, such that
all possible indices $n$ and $i\com j$ we have
\begin{equation}\label{grad}
\deg\,C^{\,(n)}=0
\quad\textrm{and}\quad
\deg\,F_{ij}^{\,(n)}=\bi\ts+\bj\,\ts.
\end{equation}
By their definition, the elements 
$C^{\ts(1)},C^{\ts(3)},\,\ldots\in\A_N$ are central.
Due to (\ref{FN-})
{\begin{equation}\label{F-}
F_{-i,-j}^{\,(n)}\ts=
(-1)^{\,n-1}\,
F_{ij}^{\,(n)}.
\end{equation}

\begin{theorem*}\label{T1.5}\ 
a) The algebra $\A_N$ is generated by the elements\/
$C^{\ts(1)},C^{\ts(3)},\,\ldots$ and\/
$F_{ij}^{\ts(1)},F_{ij}^{\ts(2)},\,\ldots\,\,$.
\\
b) The central elements 
$C^{\ts(1)},C^{\ts(3)},\,\ldots$ of\/ $\A_N$
are algebraically independent.
\\[2pt] 
c) Together with the
centrality and algebraic independence of\/
$C^{\ts(1)},C^{\ts(3)},\,\ldots\,\,$,
the defining relations of the\/ 
$\ZZ_2\ts$-\ts graded algebra\/ $\A_N$ are\/ {\rm(\ref{F-})} and
$$
[\ts F_{ij}^{\ts(m)},\ts F_{kl}^{\ts(n)}\ts]
\,=\,
F_{il}^{\ts(m+n-1)}\,\de_{kj}-
(-1)^{\,(\ts\bi\,+\ts\bj\ts)\ts(\ts\bk\,+\,\bl\ts)}\,
\de_{il}\,F_{kj}^{\ts(m+n-1)}\,\,+\vspace{2pt}
$$
$$
(-1)^{\,m\,-\ts1}\,(\,
F_{-i,l}^{\ts(m+n-1)}\,\de_{-k,j}\,-
(-1)^{\,(\ts\bi\,+\ts\bj\ts)\ts(\ts\bk\,+\,\bl\ts)}\,
\de_{i,-l}\,F_{k,-j}^{\ts(m+n-1)}\,)\,\,\,+
$$
$$
(-1)^{\ts\,\bj\,\bk\ts\,+\,\bj\,\bl\ts\,+\,\bk\,\bl}\,\,
\sum_{r=1}^{\min(m,n)-1}
(\,
F_{il}^{\ts(m+n-r-1)}\,F_{kj}^{\ts(r)}-
F_{il}^{\ts(r)}\,F_{kj}^{\ts(m+n-r-1)}
\,)\,\,\,+
$$
$$
(-1)^{\ts\,\bj\,\bk\ts\,+\,\bj\,\bl\ts\,+\,\bk\,\bl\ts\,+\,\bk\ts\,+\,\bl}
\,\,\,\times\vspace{-4pt}
$$
\begin{equation}\label{defrel}
\sum_{r=1}^{\min(m,n)-1}
(-1)^{\,m\,+\,r}\,
(\,
F_{-i,l}^{\ts(m+n-r-1)}\,F_{-k,j}^{\ts(r)}-
F_{i,-l}^{\ts(r)}\,F_{k,-j}^{\ts(m+n-r-1)}
\,)
\end{equation}
where\/ $m\com n=1\com2\com\ldots$ and\/
$i\com j\com k\com l=\pm1\lc\!\pm\!N$.
\end{theorem*}

The proof will be given in Section 3.
In particular, Theorem \ref{T1.5} shows that the algebra $\A_N$
is isomorphic to the tensor product of its two subalgebras,
generated by the elements $C^{\ts(1)},C^{\ts(3)},\,\ldots$
and by the elements $F_{ij}^{\ts(1)},F_{ij}^{\ts(2)},\,\ldots$ respectively.
Denote the latter subalgebra by $\B_N\ts$, it
is a $\ZZ_2\ts$-\ts graded subalgebra.

The algebra $\B_N$ appeared in \cite{N3} in the following guise. 
Let us consider the associative
unital $\ZZ_2$-graded algebra $\YN$ over the field $\CC$
with the countable
set of generators $T^{(n)}_{ij}$ where $n=1\com2\com\,\ldots$ and
$i,j=\pm1\lc\!\pm\!N$.
The $\ZZ_2$-gradation on the algebra $\YN$
is determined by setting $\deg\,T_{ij}^{(n)}=\bi\ts+\bj\,$ for any $n\ge1$.
To write down the defining relations for these generators, put
$$
T_{ij}(x)=
\de_{ij}\cdot1+T_{ij}^{\ts(1)}\ts x^{-1}+T_{ij}^{\ts(2)}\ts x^{-2}+\ldots
$$
where $x$ is a formal parameter, so that 
$T_{ij}(x)\in\YN\ts[[\ts x^{\ts-1}]]\ts$.
Then for all possible indices $i\com j\com k\com l$ we have the relations
\begin{equation}\label{Fu-}
T_{-i,-j}(x)=T_{ij}(-x)
\end{equation}
and
$$
(\ts x^2-y^2\ts)
\cdot
[\,T_{ij}(x)\ts,T_{kl}(y)\ts]\cdot
(-1)^{\,\bi\,\bk\ts\,+\,\bi\,\bl\ts\,+\,\bk\,\bl}\,=
\vspace{4pt}
$$
$$
(x+y)
\cdot
(\ts T_{kj}(x)\,T_{il}(y)-T_{kj}(y)\,T_{il}(x)\ts)\ -
$$
\begin{equation}\label{yangrel}
(x-y)
\cdot
(\ts T_{-k,\ts j}(x)\,T_{-i,\ts l}(y)-T_{k,\ts -j}(y)\,T_{i,\ts -l}(x) \ts)
\cdot
(-1)^{\ts\bk\,\ts+\,\bl}
\vspace{4pt}
\end{equation}
where $y$ is a formal paramater independent of $x\ts$, 
so that (\ref{yangrel}) is an equality
in the algebra of formal Laurent series 
in\/ $x^{\ts-1},y^{\ts-1}$ with coefficients in\/ $\YN\ts$.

The algebra $\YN$ is called the 
\textit{Yangian} of the Lie superalgebra $\qN\ts$.
Note that the centre of the 
associative superalgebra $\YN$ with $N\ge1$ is not trivial.
For a description of the centre of $\YN$ see
\cite[Section 3]{N3}. In particular, all central elements
of $\YN$ have $\ZZ_2$\ts-\ts degree~0.
In our Section 3 we will prove

\begin{proposition*}\label{P1.6}
The assignement\/ $F_{ij}^{\ts(n)}\mapsto(-1)^{\,\bi}\ts\,T_{ji}^{\ts(n)}$
for any\/ $n=1\com2\com\,\ldots$ and $i,j=\pm1\lc\!\pm\!N$ 
extents to an anti-isomorphism of\/ $\ZZ_2\ts$-\ts graded algebras
$$
\om:\,\B_N\ts\to\ts\YN\,.
$$
\end{proposition*}

Now denote by $\om_{N+M}$ the principal anti-automorphism of
the enveloping algebra $\UNM\ts$. It
preserves the subalgebra $\UM\subset\UNM\ts$. Hence it
also preserves the centralizer $\AMN\subset\UNM$ of that subalgebra.
For any $M\ge0$ let $\pi_M:\ts\A_N\to\AMN\ts$ be the canonical
homomorphism. By definition,
\begin{equation}\label{pim}
\pi_M\ts(\ts F_{ij}^{\ts(n)})\ts=\ts
F_{ij\ts|\ts N+M}^{\,(n)}
\end{equation}
for any $n=1\com2\com\,\ldots$ and $i\com j=\pm1\lc\!\pm\!N\ts$.
Using Proposition \ref{P1.6}, we can define a homomorphism
$\tau_M:\ts\YN\to\AMN$ by the equality 
$$
\tau_M\circ\ts\om\ts=\ts\om_{N+M}\circ\ts(\ts\pi_M\ts|\ts\B_N)\,.
$$
By (\ref{pim}),
$$
\tau_M\ts(\ts T_{ij}^{\ts(n)})
\ts=\ts
(-1)^{\ts\bj}\,\,\om_{N+M}\ts(\ts F_{ji\ts|\ts N+M}^{\ts(n)}\ts)\,.
\vspace{4pt}
$$

Using the homomorphisms $\tau_M$ for all
$M=0\com1\com2\com\,\ldots\,\ts$, one can define a family of
irreducible finite-dimensional $\YN\ts$-\ts modules\ts; see 
\cite[Section 1]{N4} and \cite{P}.
Another family of irreducible finite-dimensional 
$\YN\ts$-\ts modules 
can be defined by using the results of \cite{N1} and \cite[Section 5]{N3}.
It should be possible to give a parametrization
of all irreducible finite-dimensional $\YN\ts$-\ts modules, similarly to
the parametrization of the irreducible finite-dimensional 
${\rm Y}(\ts\mathfrak{gl}_N)\ts$-\ts modules as given by V.\,Drinfeld\ts;
see \cite[Theorem 2]{D2} and \cite{M}.

It was shown in \cite{N3} that
the associative $\ZZ_2\ts$-\ts graded algebra $\YN$
has a natural Hopf superalgebra structure. In particular, 
the homomorphism of comultiplication $\YN\to\YN\ot\YN$ can be defined by
\begin{equation}\label{comult}
T_{ij}(x)\,\mapsto\,\sum_{k}\,\,
T_{ik}(x)\ot T_{kj}(x)\cdot
{(-1)}^{\ts(\ts\bi\ts+\ts\bk\ts)(\ts\bj\ts+\ts\bk\ts)}
\end{equation}
where the tensor product is over the subalgebra 
$\CC[[\ts x^{\ts -1}]]$ of $\YN\ts[[\ts x^{\ts-1}]]\ts$, 
and $k$ runs through $\pm1\lc\!\pm\!N\ts$.
See \cite[Section 2]{N3} for the definitions of
the the counit map $\YN\to\CC$ and the antipodal map $\YN\to\YN\ts$.

There is a distinguished ascending $\ZZ\ts$-\ts filtration on the
associative algebra $\YN\ts$.
It is obtained by assigning to every generator $F_{ij}^{\ts(n)}$
the degree $n-1\ts$. 
The corresponding
$\ZZ\ts$-\ts graded algebra will be denoted by ${\rm{gr}}\,\YN\ts$.
Let $G_{ij}^{\ts(n)}$ be the element of ${\rm gr}\ts\YN$
corresponding to the generator $F_{ij}^{\ts(n)}\in\YN\ts$. 
The algebra ${\rm gr}\,\YN$ inherits the $\ZZ_2\ts$-\ts gradation 
from the algebra $\YN\ts$, so that 
$$
\deg\ts G_{ij}^{\ts(n)}=\bi\ts+\ts\bj\,.
$$
Moreover, ${\rm gr}\,\YN$ inherits from $\YN$
the Hopf superalgebra structure. It follows from
the defintion (\ref{comult}) that with respect to the
homomorphism of comultiplication 
${\rm gr}\,\YN\to{\rm gr}\,\YN\ot{\rm gr}\,\YN\ts$, 
for any $n\ge1$ we have
$$
G_{ij}^{\ts(n)}\mapsto G_{ij}^{\ts(n)}\ot1+1\ot G_{ij}^{\ts(n)}.
$$

On the other hand, 
for arbitrary Lie superalgebra $\a\ts$, a comultiplication map
$\U(\a)\to\U(\a)\ot\U(\a)$ can be defined
for $X\in\a$ by $X\mapsto X\ot1+1\ot X$, and then
extended to a homomorphism of $\ZZ_2$\ts-\ts graded associative algebras
by using the convention (\ref{otp}). Let us now consider the
enveloping algebra $\Ug$ of 
the \textit{twisted polynomial current\/} Lie superalgebra
$$
\g=\{\,X(t)\in\glN[\ts t\ts]\,:\,\et(\ts X(t))=X(-t)\,\}\,.
$$
Here we employ the automorphism (\ref{1.1}) of the 
Lie superalgebra $\glN\ts$. As a vector space, $\g$ is spanned
by the elements
\begin{equation}\label{Eu}
E_{\ts ij}\,t^{\ts n}+E_{\ts-i,-j}\ts(-\ts t)^{\ts n} 
\end{equation}
where $n=0\com1\com2\,\ldots\,$ and $i\com j=\pm1\lc\!\pm\!N\ts$. 
Note that the $\ZZ_2\ts$-\ts degree of the element (\ref{Eu})
equals $\bi\ts+\bj\,\,$. 
The algebra $\Ug$ also has a natural $\ZZ\ts$-\ts gradation,
such that the degree of the element (\ref{Eu}) is $n\ts$. 

It turns out that $\Ug$ and
${\rm gr}\,\YN$ are isomorphic as Hopf superalgebras.
By \cite[Theorem 2.3]{N3} their isomorphism
$\Ug\to{\rm gr}\,\YN$
can be defined by mapping 
the element (\ref{Eu}) of the algebra $\Ug$ to the element 
$(-1)^{\,\bi\,+1}\,G_{ji}^{\ts(n+1)}$ of the algebra ${\rm gr}\,\YN\ts$.
Moreover, $\YN$ is a deformation
of $\Ug$ as a Hopf superalgebra; see the end of
\cite[Section 2]{N3} for details.

Let us finish this introductory section with a few remarks of a
historical nature.
Our construction of the algebra $\A_N$ follows a similar construction
for the general linear Lie algebra $\mathfrak{gl}_N$ instead of the 
queer Lie superalgebra $\qN\ts$, due to G.\,Olshanski \cite{O1,O2}.
It was him who first considered the inverse limit of 
the sequence of centralizers of $\mathfrak{gl}_M$
in the enveloping algebras $\U(\ts\mathfrak{gl}_{N+M})$
for $M=1\com2\com\ts\ldots\,\ts$.
Following a suggestion
of B.\,Feigin, he then described the inverse limit in terms
of the Yangian ${\rm Y}(\ts\mathfrak{gl}_N)$ of the Lie algebra
$\mathfrak{gl}_N$. The latter Yangian
is a deformation of the
enveloping algebra of the polynomial current Lie algebra
$\mathfrak{gl}_N[t]$ in the class of Hopf algebras \cite{D1}.

The elements (\ref{FN}) of $\UN$
were initially considered by A.\,Sergeev \cite{S1}, in
order to describe the centre of the superalgebra $\UN\ts$.
The homomorphisms  $\al_M:\AMN\to\A_N^{M-1}$ for 
$M=1\com2\com\ts\ldots\,\ts$ and the elements
$F_{ij}^{\ts(1)},F_{ij}^{\ts(2)},\,\ldots$
of the inverse limit algebra $\A_N$ were introduced by M.\,Nazarov
following \cite{O1,O2}. He then identified the 
algebra \textit{defined\/} by the relations (\ref{Fu-}) and (\ref{yangrel}),
as a deformation of the enveloping algebra $\Ug$
in the class of Hopf superalgebras \cite{N3}.
It was also explained in \cite{N3} why this deformation
should be called the Yangian of $\qN$. However, our Theorems \ref{T1.1}
and \ref{T1.5} were only conjectured by M.\,Nazarov.
The purpose of the present article is to prove these conjectures.

We hope these remarks indicate importance of the
role that G.\,Olshanski
played at various stages of our work. We are very grateful to him
for friendly advice.
This work was finished while M.\,Nazarov stayed at the
Max Planck Institute of Mathematics in Bonn. He is grateful to
the Institute for hospitality.
M.\,Nazarov has been also supported by the EC grant
MRTN-CT-2003-505078.


\section{\hskip-5pt.\hskip6pt Proof of Theorem \ref{T1.1}}

We will use basic properties of complex semisimple associative
superalgebras and their modules \cite{J}. Most of these properties
were first established in \cite{W}, in a generality greater than
we need in the present article.
We will also use the following
simple lemma. Its proof carries over almost verbatim from the ungraded case,
but we shall include the proof for the sake of completeness. 
Let $\A$ be any finite dimensional $\ZZ_2$\ts-\ts graded associative algebra
over the complex field $\CC$.
Let $G$ be a finite group of automorphisms of $\A$. The crossed product
algebra $G\ltimes\A$ is also $\ZZ_2$\ts-\ts graded: for any $g\in G$
we have $\deg g=0$ in $G\ltimes\A$.

\begin{lemma*}\label{L2.1}
Suppose the superalgebra $\A$ is semisimple.
Then the superalgebra $G\ltimes A$ is also semisimple.
\end{lemma*}

\begin{proof}
We will write $G\ltimes\A=\B$. Let $V$ be any module over the superalgebra
$\B$, and let $\rho:\B\to\End(V)$ be the corresponding homomorphism. Here
we assume that the homomorphism $\rho$ preserves $\ZZ_2$\ts-\ts gradation.
Let $U\subset V$ be any $\B$\ts-\ts submodule. 
Since $\A$ is semisimple,
we have the decomposition $V=U\oplus U'$ into a direct sum
of $\A$\ts-\ts modules for some $U'\subset V$. 
Let $P\in\End(V)$ be the projection onto $U$ along $U'$. Put
$$
S\,=\,\frac1{|G|}\,\sum_{g\in G}\,\rho(g)\ts P\ts\rho(g)^{-1}\in\End(V)\,.
$$
For any element $X\in\A$ let $X^g$ be 
its image under the automorphism $g$. Then
$$
\rho(X)\ts S
\,=\,
\frac1{|G|}\,\sum_{g\in G}\,
\rho(X)\ts\rho(g)\ts P\ts\rho(g)^{-1}
\,=\,
\frac1{|G|}\,\sum_{g\in G}\,\rho(g)\ts\rho(X^g)\ts P\ts\rho(g)^{-1}
$$
$$
\,=\,
\frac1{|G|}\,\sum_{g\in G}\,\rho(g)\ts P\ts\rho(X^g)\ts\rho(g)^{-1}
\,=\,
\frac1{|G|}\,\sum_{g\in G}\,\rho(g)\ts P\ts\rho(g)^{-1}\ts\rho(X)
\,=\,
S\ts\rho(X)\,.
$$
For any $h\in G$ we also have $\rho(h)\ts S=S\ts\rho(h)$
by the definition of $S$. Note that $S\in\End(V)$ is of
$\ZZ_2$\ts-\ts degree $0$, as well as $P$ is.
So ${\rm ker\,}S\subset V$ is a $\B$\ts-\ts submodule.

Since $U$ is a $\B\ts$-\ts submodule, we have the equalities
$
P\ts\rho(g)\ts P=\rho(g)\ts P
$
for all $g\in G$. Using the definition of $S$, these equalities 
imply that $S\ts P=P$ and $P\ts S=S$. The latter pair of equalities
guarantees that ${\rm im\,}S={\rm im\,}P=U$ and $S^{\ts2}=S$.
So $V=U\oplus\ts{\rm ker\,}S$.
Since $V$ is an arbitrary module over the superalgebra $\B$, this
superalgebra is semisimple by
{[\ts J\ts,\ts Proposition 2.4\ts]}
\qed
\end{proof}

We will also need the following general ``double centralizer theorem''.
Let $V$ be any $\ZZ_2$\ts-\ts graded complex vector space.
The associative algebra $\End(V)$ is then also $\ZZ_2$\ts-\ts graded.
Take any subalgebra $\B$ in the superalgebra $\End(V)$.
Here we assume that $\B$ as a vector space splits 
into the direct sum of its subspaces of $\ZZ_2$\ts-\ts degrees $0$ and $1$.
Denote by $\Bp$ the centralizer of $\B$ in the superalgebra $\End(V)$.

\begin{proposition*}\label{P2.2}
Suppose that the superalgebra $\B$ is finite dimensional and semisimple. 
Also suppose that the $\Bp$\ts-\ts module $V$ is finitely generated.
Then $\B=\Bpp$ in $\End(V)$.
\end{proposition*}

\begin{proof}
We shall prove that for any homogeneous $v_1\lc v_n\in V$
and $X\in\Bpp$, there exists $Y\in\B$ such that $Xv_r=Yv_r$ for any index
$r=1\lc n$. Then we will choose $v_1\lc v_n$ to be homogeneous generators
of $V$ over $\Bp$. By writing any vector $v\in V$ as the sum
$Z_1v_1+\ldots+Z_nv_n$ for some homogeneous $Z_1\lc Z_n\in\Bp$,
we will get
\begin{eqnarray*}
Xv\,&=&\,\sum_{r=1}^n\,XZ_r\ts v_r
\,=\,
\sum_{r=1}^n\,\,
(-1)^{\,\deg X\,\deg Z_r}\,
Z_r\ts X\ts v_r
\\
&=&\,
\sum_{r=1}^n\,\,
(-1)^{\,\deg X\,\deg Z_r}\,
Z_r\ts Y\ts v_r
\,=\,
\sum_{r=1}^n\,YZ_r\ts v_r\,=\,Yv\,.
\end{eqnarray*}
\nopagebreak
Along with the obvious embedding $\B\subset\Bpp$,
this will prove Proposition~\ref{P2.2}.

Recall that the $\ZZ_2$\ts-\ts graded vector space $\Vb$ 
{\it opposite\/} to $V$ is obtained from $V$ by changing the
$\ZZ_2$\ts-\ts gradation deg to \,deg\ts+1. Define the action 
of the algebra $\B$ in $\Vb$ as the pullback its action in $V$
via the involutive automorphism $Y\mapsto(-1)^{\deg Y}\,Y$,
where $Y$ is any homogeneous element of the $\ZZ_2$\ts-\ts graded 
algebra $\B$. Now consider the direct sum of $\B$\ts-\ts modules
$$
W\,=\,\bigoplus_{r=1}^{\ts n}\,V_r\,,
$$ 
where the $\B$\ts-\ts module $V_r$ equals $V$ or $\Vb$
depending on whether $\deg v_r$ in $V$ is~$0$~or~$1$.
For each index $r=1\lc n$ we have an embedding of vector spaces
$A_r:V\to W$, and a projection $B_r:W\to V$.
Note that the $\ZZ_2$\ts-\ts degrees
of the linear maps $A_r$ and $B_r$ coincide with that of the vector $v_r$.
We also have the equality 
$B_p\ts A_q=\de_{pq}\cdot\id$ in $\End(V)$,
and the equality
$$
\sum_{r=1}^n\,A_r\ts B_r\,=\,\id
$$
in $\End(W)$\,. Any homogeneous element $Y\in\B$ acts in $W$ as
the linear operator
$$
\Yt\,=\,\sum_{r=1}^n\,\,
(-1)^{\,\deg v_r\,\deg Y}\,
A_r\ts Y\ts B_r\,.
$$ 
Given $X\in\Bpp$, put
$$
\Xt\,=\,\sum_{r=1}^n\,\,
(-1)^{\,\deg v_r\,\deg X}\,
A_r\ts X\ts B_r\,.
$$

Any homogeneous element $Z\in\End(W)$ can be written as
$$
Z\,=\sum_{p,q=1}^n\,A_p\ts Z_{pq}\ts B_q
$$
for some homogeneous $Z_{pq}\in\End(V)$, where
$\deg Z_{pq}=\deg Z+\deg v_p+\deg v_q$. Then
\begin{eqnarray}
\Yt Z&=&
\sum_{r=1}^n\,\,
(-1)^{\,\deg v_r\,\deg Y}\,
A_r\ts Y\ts B_r\,
\cdot\sum_{p,q=1}^n A_p\ts Z_{pq}\ts B_q
\label{2.1}
\\
&=&
\sum_{p,q=1}^n
(-1)^{\,\deg v_p\,\deg Y}\,
A_p\ts Y Z_{pq}\ts B_q\,,
\nonumber
\\
Z\,\Yt&=&
\sum_{p,q=1}^n A_p\ts Z_{pq}\ts B_q
\,\cdot\,
\sum_{r=1}^n\,\,
(-1)^{\,\deg v_r\,\deg Y}\,
A_r\ts Y\ts B_r
\label{2.2}
\\
&=&
\sum_{p,q=1}^n
(-1)^{\,\deg v_q\,\deg Y}\,
A_p\ts Z_{pq}\ts Y\ts B_q\,.
\nonumber
\end{eqnarray}
Suppose that the element $Z\in\End(W)$ belongs to the centralizer
of the image of the superalgebra $\B$ in $\End(W)$.
Due to (\ref{2.1}) and (\ref{2.2}), the assumption $[\ts\Yt,Z\ts]=0$ is 
equivalent
to the collection of equalities $[\ts Y, Z_{pq}\ts]=0$
for all $p\com q=1\lc n$.
Since here $Y\in\B$ is arbitrary and $X\in\Bpp$,
we then have $[\ts X\com Z_{pq}\ts]=0$
for all $p\com q=1\lc n$. A calculation similar to
(\ref{2.1}) and (\ref{2.2}) then shows that $[\ts\Xt\com Z\ts]=0$ in 
$\End(W)$.

Now take the vector
$$
w\,=\,\sum_{r=1}^n\,A_r\ts v_r\in W\,,
$$
it has $\ZZ_2$\ts-\ts degree $0$. Let $U$ be the cyclic span of
the vector $w$ under the action of $\B$. Since the superalgebra $\B$
is finite dimensional semisimple, we have the decomposition
$W=U\oplus U'$ into direct sum of $\B$\ts-\ts modules for some $U'\subset W$;
see [\ts J\ts,\ts Proposition~2.4].
Choose the element $Z\in\End(W)$ from the centralizer of the image of $\B$,
to be the projector onto $U'$ along $U$. Here
$\deg Z=0$. Then $Z\Xt w=\Xt Zw=0$, and $\Xt w\in U$.
So there exists $Y\in\B$ such that $\Xt w=\Yt w$.
The last equality means that $X v_r=Y v_r$ for each $r=1\lc n$.
\qed
\end{proof}

In the notation of Proposition \ref{P2.2}, we have the following corollary.

\begin{corollary*}\label{C2.3}
Suppose that the vector space $V$ is finite dimensional,
and that the superalgebra $\B\subset\End(V)$ is semisimple.
Then $\B=\Bpp$ in\/ $\End(V)$.
\end{corollary*}

Now for any integer $n\ge1$ consider the tensor product 
$(\ts\End(V))^{\ot\ts n}$ of $n$ copies of the
$\ZZ_2$\ts-\ts graded algebra $\End(V)$.
This tensor product is a $\ZZ_2$\ts-\ts graded
associative algebra defined using the conventions
(\ref{otp}) and (\ref{otd}).
The proof of the following lemma is also included
for the sake of completeness.

\begin{lemma*}\label{L2.0}
Suppose that $\B$ contains the identity $1\in\End(V)$.
Then the centralizer of\/ $\B^{\ot\ts n}$
in the superalgebra $(\ts\End(V))^{\ot\ts n}$ coincides with
$(\ts\Bp)^{\ot\ts n}$.
\end{lemma*}

\begin{proof}
We will use the induction on $n\ts$. If $n=1$, the statement of
Lemma \ref{L2.0} is tautological. Suppose that $n>1$.
The centralizer of $\B^{\ot\ts n}$
contains $(\ts\Bp)^{\ot\ts n}$
due to the conventions (\ref{otp}) and (\ref{otd}).
Now suppose that for some homogeneous elements
$X_1\lc X_l\in\End(V)$
and
$Y_1\lc Y_l\in(\ts\End(V))^{\ot\ts(n-1)}$
the sum
$X_1\ts Y_1+\ldots+X_l\ts Y_l$
belongs to the centralizer of $\B^{\ot\ts n}$. 
In particular, then for any homogeneous $Y\in\B^{\ot\ts(n-1)}$ we have
$$
\begin{aligned}[t]
\sum_{k=1}^l\,X_k\ot(\ts Y_k\ts Y)
&\,=\,
\sum_{k=1}^l\,(\ts X_k\ot Y_k)\ts(\ts1\ot Y)
\\
&\,=\,
\sum_{k=1}^l\,
(-1)^{\,(\deg X_k\ts+\,\deg Y_k)\,\deg Y}
(\ts1\ot Y)\ts(\ts X_k\ot Y_k)
\\
&\,=\,
\sum_{k=1}^l\,
(-1)^{\,\deg Y_k\,\deg Y}
X_k\ot(Y\ts Y_k)\,.
\end{aligned}
$$
We may assume that the elements $X_1\lc X_l$ are linearly
independent, then the above equalities imply that for any
$k$ the element $Y_k\in(\ts\End(V))^{\ot\ts(n-1)}$
belongs to the centralizer of $\B^{\ot\ts (n-1)}$. 
Then $Y_k\in(\ts\Bp)^{\ot\ts (n-1)}$ by the induction assumption.
A similar argument shows that $X_k\in\Bp$ for any index $k$.
\qed
\end{proof}

For any $N\ge0$ and $M\ge1$ take the $\ZZ_2$\ts-\ts graded
vector space $\CNM$. Identify the $\ZZ_2$\ts-\ts 
graded vector spaces $\CN$ and $\CM$ with the subspaces in 
$\CNM$ spanned by the vectors $e_i$
where $i=-N\lc\!-\!1\com1\lc N$ and 
$i=-N\!-\!M\lc\!-\!N\!-1\com N+1\lc N+M$ respectively.
The decomposition 
\begin{equation}\label{dec}
\CN=\,\CN\op\CM
\end{equation}
determines the embeddings of the Lie superalgebras 
$\glN$ and $\glM$ into $\glNM$, and of their subalgebras
$\qN$ and $\qM$ into the Lie algebra $\qNM\ts$.

Let us now regard $\End(\CN)$ and $\Q_M$ as subalgebras in the associative
superalgebra $\End(\CNM)\ts$, using the
decomposition (\ref{dec}). The elements
$$
\sum_{|i|>N}E_{ii}
\quad\text{and}\quad
\sum_{|i|>N}(-1)^{\ts\bi}\,E_{i,-i}
$$
of $\End(\CNM)$ span a subalgebra, isomorphic to 
the associative algebra $\Q_1$. Using this isomorphism,
a direct calculation shows that the centralizer of 
$\Q_M$ in $\End(\CNM)$ coincides with $\End(\CN)\op\Q_1\ts$.
This centralizer is a semisimple associative
$\ZZ_2$\ts-\ts graded algebra, denote it by $\C$. 

For any integer $n\ge1$ consider the tensor product
$V=(\CNM)^{\ts\ot\ts n}$ of $n$ copies
of the $\ZZ_2$\ts-\ts graded vector space $\CNM$.
Identify the algebras
\begin{equation}\label{ide}
(\ts\End(\CNM))^{\ot\ts n}
\quad\text{and}\quad
\End((\CNM)^{\ot\ts n})
\end{equation}
so that for any homogeneous elements $X_1\lc X_n\in\End(\CNM)$
and any homogeneous vectors $u_1\lc u_n\in\CNM$
$$
(\ts X_1\ot\ldots\ot X_n)\,(\ts u_1\ot\ldots\ot u_n)=
(-1)^d\,X_1 u_1\ot\ldots\ot X_n u_n
$$
where 
$$
d\ =\!\sum_{1\le p<q\le n}\!\deg\ts u_p\,\,\deg\ts X_{q}\,.
$$
By identifying the two algebras we
determine an action on $V$ of the subalgebra
\begin{equation}\label{suc}
\C^{\ts\ot\ts n}\subset(\ts\End(\CNM))^{\ot\ts n}.
\end{equation}

The symmetric group $S_n$ acts on $V$ so that
for any adjacent transposition $\sigma_p=(p\ts,p+1)$ 
$$
\sigma_p\,(\ts u_1\ot\ldots\ot u_p\ot u_{p+1}\ot\ldots\ot u_n)=
\hskip-4pt
$$
$$
(-1)^{\,\deg u_p\,\deg u_{p+1}}\,
u_1\ot\ldots\ot u_{p+1}\ot u_p\ot\ldots\ot u_n\,.
$$
The action of $S_n$ and $\C^{\ts\ot\ts n}$ on $V$
extends to that of $S_n\ltimes\C^{\ts\ot\ts n}$.
Here the group $S_n$ acts by automorphisms of the algebra 
$(\ts\End(\CNM))^{\ot\ts n}$ so that
$$
\sigma_p\,(\ts X_1\ot\ldots\ot X_p\ot X_{p+1}\ot\ldots\ot X_n)=
$$
\begin{equation}\label{six}
(-1)^{\,\deg X_p\,\deg X_{p+1}}\,
X_1\ot\ldots\ot X_{p+1}\ot X_p\ot\ldots\ot X_n\,.
\end{equation}

Now consider the action of enveloping algebra $\U(\qN)$ on 
the vector space $V$, 
as of a subalgebra of the $\ZZ_2$\ts-\ts graded 
associative algebra $\U(\glNM)\ts$. We use
the comultiplication 
\begin{equation}\label{com}
\U(\glNM)\to\U(\glNM)^{\ot\ts n}
\end{equation}
along with the identification of the two algebras in (\ref{ide}).

\begin{proposition*}\label{P2.4}
The centralizer of the image of\/ $\U(\qM)$ in\/
the associative superalgebra\/
$\End(V)$ coincides with the image of\/ $S_n\ltimes\C^{\ts\ot\ts n}$.
\end{proposition*}

\begin{proof}
The $\ZZ_2$\ts-\ts graded algebra $\C^{\ts\ot\ts n}$ is 
semisimple, see \cite[Proposition 2.10]{J}. So is
the crossed product $S_n\ltimes\C^{\ts\ot\ts n}$, see
our Lemma \ref{L2.1}. Let $\B$ be the image of the crossed product 
in $\End(V)$.
The $\ZZ_2$\ts-\ts graded algebra $\B$ is semisimple too.
By Corollary~\ref{C2.3} it suffices
to prove that the centralizer of $\B$ in $\End(V)$
coincides with the image of $\U(\qM)\ts$.

Let us describe the latter image. Consider the span in
$\End(\CNM)$ of $\qM$ and of the identity element $1$. It is
a subalgebra in $\End(\CNM)\ts$, denote this subalgebra by $\D$.
We will prove that the invariant subalgebra
\begin{equation*}
(\ts\D^{\ot\ts n})^{S_n}\subset(\ts\End(\CNM))^{\ot\ts n}
\end{equation*}
coincides with the image of $\U(\qM)\ts$.
Here we use the definition (\ref{six}), and
the comultiplication (\ref{com}). There is no need
to identify the two algebras (\ref{ide}) here.
 
Due to the convention (\ref{otp}) and the
definition of the comultiplication (\ref{com}),
the image of $\U(\qM)$ in $(\ts\End(\CNM))^{\ot\ts n}$
is contained in $(\ts\D^{\ot\ts n})^{S_n}$. Now for any $n$ elements
$X_1\lc X_n\in\D$ 
consider their \textit{symmetrized\/} tensor product 
\begin{equation}\label{ss}
\langle\ts X_1\lc X_n\rangle=
\frac{1}{n!}\sum_{\sigma\in S_n}
\sigma\,(\ts X_1\ot\ldots\ot X_n)\,.
\end{equation}
Suppose that for some $p\in\{1\lc n\}$ we have $X_q\in\CC\,1$
if and only if $p<q\ts$. By induction on $p\ts$, let us prove that
$
\langle\ts X_1\lc X_p,1\lc1\rangle
$
belongs to the image of $\U(\qM)$ in $(\ts\End(\CNM))^{\ot\ts n}$.
This is evident if $p=1$. If $p>1$, then
$$
\langle\ts X_1\lc X_{p-1},1\lc1\rangle
\,
\langle\ts X_p,1\lc1\rangle
$$
equals
$$
\frac{n-p+1}n\,\langle\ts X_1\lc X_p,1\lc1\rangle
$$
plus a sum of certain symmetrized tensor products in 
$(\ts\End(\CNM))^{\ot\ts n}$
which belong to the image of $\U(\qM)$ by the induction assumption.

To prove that the centralizer of $\B$ in $\End(V)$
coincides with the image of $\U(\qM)\ts$, it now suffices
to show that the centralizer of the subalgebra (\ref{suc})
coincides with $\D^{\ot\ts n}$. But we have 
$\C^{\ts\prime}=\D$ in $\End(\CNM)$ by definition. 
Using Lemma \ref{L2.0} we now complete the proof
of Proposition 2.5.
\qed
\end{proof}

Consider the symmetric algebra $\S(\qNM)$
of the Lie superalgebra $\qNM\ts$. The standard filtration 
\begin{equation}\label{filtr}
\CC=\U^{\ts0}(\qNM)
\subset
\U^{\ts 1}(\qNM)
\subset
\U^{\ts 2}(\qNM)
\subset
\ldots
\end{equation}
on the algebra $\UNM$ determines for each $n=1,2,\ldots$
a linear map
\begin{equation}\label{ison}
\U^{\ts n}(\qNM)\,/\,\U^{\ts n-1}(\qNM)\,\to\,\S^{\ts n}(\qNM)\,,
\end{equation}
which is bijective
by the Poincar\'e\,-Birkhoff\,-Witt theorem for Lie \text{superalgebras} 
\cite[Theorem 5.15]{MM}. Using (\ref{ison}), for all
indices $i,j=\pm1\lc\!\pm\!(N+M)$
define an element $f_{ij\ts|\ts N+M}^{\,(n)}\in\S^{\ts n}(\qNM)$ 
as the image of the element
$$
F_{ij\ts|\ts N+M}^{\,(n)}\in\U^{\ts n}(\qNM)\,;
$$
the latter element is the sum (\ref{FN}) where
each of the indices $k_1\lc k_{n-1}$~runs through $\pm1\lc\!\pm\!(N+M)\ts$.
Define $c_{\ts N+M}^{\,(n)}\in\S^{\ts n}(\qNM)$  
as the image of
$$
C_{\ts N+M}^{\,(n)}\in\U^{\ts n}(\qNM)\,;
$$
see (\ref{CN}).
Note that if $n$ is even, then
$C_{\ts N+M}^{\,(n)}=0$ and hence $c_{\ts N+M}^{\,(n)}=0$.

Now consider adjoint action of the Lie superalgebra $\qNM$
on $\S(\qNM)\ts$. 
In particular, consider the action of $\qM$
on $\S(\qNM)$ as that of a subalgebra of $\qNM\ts$. Then 
take the subalgebra of invariants $\S(\qNM)^{\ts\qM}\subset\S(\qNM)\ts$. 

\begin{proposition*}\label{P2.5}
The subalgebra\/ $\S(\qNM)^{\ts\qM}$ 
is generated by the elements 
$c_{\ts N+M}^{\,(1)}\,,\,c_{\ts N+M}^{\,(3)}\,,\,\ldots$
and
$f_{\ts ij\ts|\ts N+M}^{\,(1)}\,,\,f_{\ts ij\ts|\ts N+M}^{\,(2)}\,,\,\ldots\,$
where
$|i|\com|j|\le N$.
\end{proposition*}

\begin{proof}
Consider $(\ts\End(\CNM))^{\ot\ts n}$ as a vector space.
Define a linear map $\ph_n$ from this space to the
$n\ts$th symmetric power $\S^{\ts n}(\qNM)$ by setting
\begin{equation}\label{phin}
\ph_n\ts(\ts E_{i_1j_1}\ot\ldots\ot E_{i_nj_n})=
F_{i_1j_1}\ldots F_{i_nj_n}
\end{equation}
for any indices $i,j=\pm1\lc\!\pm\!(N+M)\ts$. The
map $\ph_n$ commutes with the natural action 
of the Lie superalgebra $\qNM$ on
$(\ts\End(\CNM))^{\ot\ts n}$ and $\S^{\ts n}(\qNM)\ts$.
The map $\ph_n$ has a right inverse linear map $\ps_n$
which commutes with the action of
$\qNM$ as well. Namely, using (\ref{ss}) we set
$$
\ps_n\ts(\ts F_{i_1j_1}\ldots F_{i_nj_n})=2^{-n}\,
\langle\ts F_{i_1j_1}\lc F_{i_nj_n}\rangle\,.
$$
It follows that the subspace of $\qM\ts$-\ts invariants
of the $\ZZ\ts$-\ts degree $n$ in $\S(\qNM)\ts$,
$$
\S^{\ts n}(\qNM)^{\ts\qM}=
\ph_n\ts(((\ts\End(\CNM))^{\ot\ts n})^{\ts\qM})=
\ph_n(\ts\B)
$$
by Proposition \ref{P2.4}\ts; here 
$\B$ is the image in of the crossed product 
$S_n\ltimes\C^{\ts\ot\ts n}$~in
$$
(\ts\End(\CNM))^{\ot\ts n}
=
\End((\CNM)^{\ot\ts n})\,.
$$

The vector subspace $\C\subset\End(\CNM)$
is spanned by the identity element $1$, 
the elements $E_{ij}$ where $|i|\com|j|\le N$
and by the element
$$
J=\sum_i\,(-1)^{\ts\bi}\,E_{i,-i}
$$
where the summation index $i$ runs through $\pm1\lc\!\pm\!(N+M)\ts$.
For each $p=1\lc n$ introduce the elements of the algebra 
$(\ts\End(\CNM))^{\ot\ts n}$,
$$
E^{\ts(p)}_{ij}=1^{\ot\ts(p-1)}\ot E_{ij}\ot1^{\ot\ts(n-p)}
\quad\text{and}\quad
J_p=1^{\ot\ts(p-1)}\ot J\ot1^{\ot\ts(n-p)}\,.
$$
The vector subspace $\B\subset(\ts\End(\CNM))^{\ot\ts n}$
is spanned by the products of the form
\begin{equation}\label{JEH}
J_{q_1}\ldots\,J_{q_b}\,
E^{\ts(p_1)}_{i_1j_1}\ldots\,E^{\ts(p_a)}_{i_aj_a}\,H
\end{equation}
where 
$$
1\le p_1<\ldots<p_a\le n
\ts,\quad
1\le q_1<\ldots<q_b\le n
\ts,
$$
$$
\{p_1\lc p_a\}\cap\{q_1\lc q_b\}=\varnothing
\ts,
$$
$$
|i_1|\com|j_1|\lc|i_a|\com|j_a|\le N
$$
and $H$ is the image in $(\ts\End(\CNM))^{\ot\ts n}$
of some permutation from $S_n\ts$.

By the definintion of the symmetric algebra $\S(\qNM)$ and
due to (\ref{six}), we have the identity 
$\ph_n\circ\ts\si=\ph_n$ for any $\si\in S_n\ts$.
Therefore it suffices to compute the $\ph_n\ts$-\ts image
of the element (\ref{JEH}) where the factor $H$ corresponds 
to a permutation of the form
$$
(1,2\lc r_1)\,(r_1+1,r_1+2\lc r_2)\,\ldots\ts(r_c+1,r_c+2\lc n)
$$
where $c\ge0$ and $0<r_1<r_2<\ldots<r_c<n\ts$. But 
for any $p=1\lc n-1$ we have the relation
\begin{equation}\label{pnp}
\ph_p(\ts X)\,\ts\ph_{n-p}(\ts Y)=\ph_n(\ts X\ot Y)\,,
\end{equation}
$$
X\in(\ts\End(\CNM))^{\ot\ts p}
\quad\textrm{and}\quad
Y\in(\ts\End(\CNM))^{\ot\ts(n-p)}\,.
$$
Hence it suffices to consider only the case where $H$
corresponds to the cyclic permutation $(1,2\lc n)\ts$.
Suppose this is the case.
Then we may assume that $p_1=1$ or $a=0\ts$.
Here we use the identity $\ph_n\circ\ts\si=\ph_n$ for 
$\si=(1,2\lc n)\ts$.

For the cyclic permutation $(1,2\lc n)$ we have
\begin{equation}\label{H}
H\,=\sum_{k_1,...,k_n}
(-1)^{\,\bk_1\,+\,\ldots\,+\,\bk_{n-1}}\,
E_{\ts k_nk_1}^{\,\ts(1)}
E_{\ts k_1k_2}^{\,\ts(2)}
\ldots
E_{\ts k_{n-1}k_n}^{\,\ts(n)}
\end{equation}
where each of the indices $k_1\lc k_n$ runs through $\pm1\lc\pm(N+M)\ts$.
Hence
\begin{equation}\label{EH}
E^{\ts(1)}_{ij}
H\,=
\sum_{k_1,...,k_{n-1}}\,
(-1)^{\,\bk_1\,+\,\ldots\,+\,\bk_{n-1}}
E_{\ts k_nk_1}^{\,\ts(1)}
E_{\ts k_1k_2}^{\,\ts(2)}
\ldots
E_{\ts k_{n-1}k_n}^{\,\ts(n)}\ts.
\end{equation}
Further, for any $p=1\lc n-1$ we have the equality
$$
E^{\ts(1)}_{ij}E^{\ts(p+1)}_{kl}H=
(-1)^{\,\bj\,\bk\,+\,\bj\,\bl\,+\,\bk\,\bl}\,
E^{\ts(1)}_{il}E^{\ts(p+1)}_{kj}H'\ts H''
$$
where the factors $H'$ and $H''$ are the images
in $(\ts\End(\CNM))^{\ot\ts n}$ of the cyclic permutations
$(1\lc p)$ and $(p+1\lc n)\ts$. Using this equality together
with the relation
(\ref{pnp}), we reduce the case $a>1$ to the case $a=1$.

Suppose that $a=1$ and $p_1=1$. In this case,
the $\ph_n\ts$-\ts image of (\ref{JEH}) equals
$$
\sum_{k_1,...,k_{n-1}}
(-1)^{\,\bk_1\,+\,\ldots\,+\,\bk_{n-1}
\,+\,
b\,\bi\,+\,q_1\,+\,\ldots\,+\,q_b}\,
F_{\ts \ts(-1)^b\ts i,k_1}
F_{\ts k_1k_2}
\ldots
F_{\ts k_{n-1}j}
$$
\vskip-12pt
$$
=(-1)^{\,b\,\bi\,+\,q_1\,+\,\ldots\,+\,q_b}\,
f_{\ts(-1)^b\ts i\ts,\ts j\ts|\ts N+M}^{\,(n)}\,;
$$
we use the equality (\ref{EH}), the definition (\ref{phin}) and
the identity $F_{-k,-l}=F_{kl}\ts$. 

It remains to consider the case $a=0$.
Then the $\ph_n\ts$-\ts image of (\ref{JEH}) equals
$$
(-1)^{\,q_1\,+\,\ldots\,+\,q_b}
\sum_{k_1,...,k_n}
(-1)^{\,\bk_1\,+\,\ldots\,+\,\bk_{n-1}\,+\,b\,\bk_n}
\,
F_{\ts(-1)^b\ts k_n,k_1}
F_{\ts k_1k_2}
\ldots
F_{\ts k_{n-1}k_n}\,;
$$
we use the equality (\ref{H}), the definition (\ref{phin}) and
the identity $F_{-k,-l}=F_{kl}\ts$. Denote by $f$ the sum
over $k_1\lc k_n$ in the above display.
If the number $b$ is even, then 
$
f=c_{\ts N+M}^{\,(n)}\,.
$
Now suppose that the number $b$ is odd, so that $f$ is
\begin{equation}\label{fbo}
\sum_{k_1,...,k_n}
(-1)^{\,\bk_1\,+\,\ldots\,+\,\bk_{n-1}\,+\,\bk_n}
\,
F_{\ts-k_n,k_1}
F_{\ts k_1k_2}
\ldots
F_{\ts k_{n-1}k_n}\,.
\end{equation}
By changing the signs of the indices $k_1\lc k_n$ in (\ref{fbo})
and then using the identity $F_{-k,-l}=F_{kl}$ one shows that
$f=(-1)^{\ts n}f\ts$. So $f=0$ if $n$ is odd. Since 
$$
\deg\, F_{\ts-k_n,k_1}=\bk_1+\bk_n+1=
\deg\ts(\ts F_{\ts k_1k_2}
\ldots
F_{\ts k_{n-1}k_n})+1\ts,
$$
the element (\ref{fbo}) of the symmetric algebra $\S(\qNM)$
can be also written as
\begin{equation}\label{fbon}
\sum_{k_1,...,k_n}
(-1)^{\,\bk_1\,+\,\ldots\,+\,\bk_{n-1}\,+\,\bk_n}
\,
F_{\ts k_1k_2}
\ldots
F_{\ts k_{n-1}k_n}\,
F_{\ts-k_n,k_1}\,.
\end{equation}
By changing the signs of the indices $k_2\lc k_n$ in (\ref{fbon})
and then using the identity $F_{-k,-l}=F_{kl}$ one shows that
$f=(-1)^{\ts n-1}f\ts$. So $f=0$ if $n$ is even.~\qed
\end{proof}

By using the fact that for each $n\ge1$
the linear map (\ref{ison}) commutes with
the adjoint action of the Lie superalgebra $\qNM$ on 
$\U^{\ts n}(\qNM)$ and $\S^{\ts n}(\qNM)\ts$,
we can now complete the proof of Theorem \ref{T1.1}.
Take any element $X\in\AMN\ts$. 
We have $X\in\U^{\ts n}(\qNM)$ for some $n\ge0$.
Let us demonstrate by induction on $n$ that $X$
belongs to the subalgebra in $\UNM$ generated by the elements
$C_{N+M}^{\,(1)}\,,C_{N+M}^{\,(3)}\,,\,\ldots$
and {\rm(\ref{elements})}. Since $\U^{\ts 0}(\qNM)=\CC$,
here we can take $n\ge1$ and make the induction assumption.
By Proposition \ref{P2.5}, the image of $X$ in 
$\S^{\ts n}(\qNM)$ under the map (\ref{ison})
is a linear combination of the products
of the elements
$c_{\ts N+M}^{\,(1)}\,,\,c_{\ts N+M}^{\,(3)}\,,\,\ldots$ 
and
$f_{\ts ij\ts|\ts N+M}^{\,(1)}\,,\,
f_{\ts ij\ts|\ts N+M}^{\,(2)}\,,\,\ldots\,$
where
$|i|\com|j|\le N$.
By replacing these elements of $\S(\qNM)$ 
respectively by
$C_{N+M}^{\,(1)}\,,C_{N+M}^{\,(3)}\,,\,\ldots$
and {\rm(\ref{elements})} in the linear combination,
we obtain a certain element $Y\in\UNM$
such that $X-Y\in\U^{\ts(n-1)}(\qNM)\ts$.
We also have $Y\in\AMN\ts$.
By applying the induction assumption to the difference
$X-Y$, we complete the the proof.


\section{\hskip-5pt.\hskip6pt Proof of Theorem \ref{T1.5}} 

Here we prove Theorem \ref{T1.5} along with Lemma \ref{L1.3}
and Propositions \ref{P1.4}\,,\ts\ref{P1.6}.

\medskip\noindent
\textit{Proof of Lemma 1.3.\/}
Together with the right ideal $\I_{\ts N+M}$ generated by
the elements (\ref{JM}), consider
the left ideal $\J_{\ts N+M}$ in $\UNM$
generated by the elements
$$
F_{\pm1,\ts N+M}
\ts\lc
F_{\pm(N+M),\ts N+M}
\,\ts.
$$
By the Poincar\'e\,-Birkhoff\,-Witt theorem 
for Lie superalgebras, 
every element $X\in\UNM$ can be uniquely written as a sum
of the products of the form
$$
F_{N+M,\ts1}^{\,p_{\ts1}}\,
F_{N+M,\ts-1}^{\,p_{\ts-1}}
\,\ldots\,
F_{N+M,\ts N+M-1}^{\,p_{\ts N+M-1}}\,
F_{N+M,\ts-N-M\ts+\ts1}^{\,p_{\ts-N-M\ts+\ts1}}
\ \,\times
$$
$$
F_{N+M,\ts N+M}^{\,p_{\ts N+M}}\,\ts
F_{N+M,\ts -N-M}^{\,p_{\ts-N-M}}\,\ts
Y\ \,\times
$$
\begin{equation}\label{HC}
F_{1,\ts N+M}^{\,q_{\ts1}}\,
F_{-1,\ts N+M}^{\,q_{\ts-1}}
\,\ldots\,
F_{N+M-1,\ts N+M}^{\,q_{\ts N+M\ts-\ts1}}\,
F_{-N-M\ts+\ts1,\ts N+M}^{\,q_{\ts-N-M\ts+\ts1}}
\end{equation}
where each of the exponents $p_k$ and $q_k$ runs through
$0\com1\com2\com\,\ldots$ or through $0\com1$ if
$k>0$ or $k<0$ respectively, whereas the factor 
$Y\in\U(\q_{\ts N+M-1})$ depends on these exponents.
Note that here
\begin{equation}\label{together}
[\,F_{N+M,\ts N+M}\ts\com\ts Y\ts]\,=\,0
\quad\textrm{and}\quad
[\,F_{-N-M,\ts N+M}\ts\com\ts Y\ts]\,=\,0\,.
\end{equation}

Now suppose that $X\in\AMN\ts$.
In particular, then we have
\begin{equation}\label{condition}
[\,F_{N+M,\ts N+M}\ts\com\ts X\ts]\,=\,0
\end{equation}
since $M\ge1$ by our assumption.
Note that if $|\ts k\ts |<N+M\ts$, then due to (\ref{FR})
$$
\begin{aligned}[t]
&
[\,F_{N+M,\ts N+M}\ts\com\ts F_{N+M,\ts k}\ts]
\,=\,F_{N+M,k}\ ,
\\[3pt]
&
[\,F_{N+M,\ts N+M}\ts\com\ts F_{k,\ts N+M}\ts]
\,=\,-\,
F_{N+M,k}\ .
\end{aligned}
$$
If $|\ts k\ts |=N+M\ts$, then
$$
[\,F_{N+M,\ts N+M}\ts\com\ts F_{k,\ts N+M}\ts]\,=\,0\,.
$$
Hence the condition (\ref{condition}) implies that $X$
is a sum of the products (\ref{HC}) where
$$
p_{\ts1}+p_{\ts-1}+\ldots+p_{\ts N+M-1}+p_{\ts-N-M\ts+\ts1}=
$$

\vspace{-16pt}
\begin{equation}\label{pq}
q_{\ts1}+q_{\ts-1}+\ldots+q_{\ts N+M-1}+q_{\ts-N-M\ts+\ts1}\ .
\end{equation}

The intersection $\I_{\ts N+M}\cap\AMN$ consists of those
elements $X\in\AMN$ which are sums of the products (\ref{HC}) where 
$$
p_{\ts1}+p_{\ts-1}+\ldots+
p_{\ts N+M}+p_{\ts-N-M}>0\,.
$$
Due to the equality (\ref{pq}), the latter inequality is equivalent to
$$
q_{\ts1}+q_{\ts-1}+\ldots+q_{\ts N+M-1}+q_{\ts-N-M\ts+\ts1}+
p_{\ts N+M}+p_{\ts-N-M}>0\,.
$$
So by using (\ref{together}), 
\begin{equation}\label{intersect}
\I_{\ts N+M}\cap\AMN=\J_{\ts N+M}\cap\AMN\,.
\end{equation}
In particular, the intersection $\I_{\ts N+M}\cap\AMN$
is a two\ts-\ts sided ideal of $\AMN\ts$. Thus we
get Part (a) of Lemma \ref{L1.3}. 

Furthemore, due to (\ref{pq}) there is only one
summand (\ref{HC}) of $X\in\AMN$ with
\begin{equation}\label{summand}
p_{\ts1}+p_{\ts-1}+\ldots+
p_{\ts N+M}+p_{\ts-N-M}=0\,,
\end{equation}
this summand has the form of $Y\in\U(\q_{\ts N+M-1})\,$.
Note that the right ideal $\I_{\ts N+M}$ of $\ts\UNM$ 
is stable under the adjoint action of the subalgebra
$\q_{\ts N+M-1}\subset\qNM\,$. Indeed, if
$|i|,|j|<N+M$ then by (\ref{FR}) we have
$$
[\,F_{ij}\ts\com\ts F_{N+M,\ts l}\,]
\,=\,-\,
(-1)^{\,\ts\bi\,\ts\bl\,+\ts\bj\,\ts\bl}\,\,\de_{\ts il}\,F_{N+M,\ts j}\,-\,
(-1)^{\,\ts\bi\,\ts\bl\,+\ts\bj\,\ts\bl}\,\,\de_{\ts i,-l}\,F_{N+M,\ts-j}
$$
for any index $l=\pm1\lc\!\pm(N+M)\ts$. 
In particular, $\I_{\ts N+M}$ is stable
under the adjoint action of $\q_{\ts M-1}\,$.
So the condition $[\ts Z\com\ts X\ts]=0$ on $X$
for any $Z\in\q_{\ts M-1}$
implies the condition $[\ts Z\com\ts Y\ts]=0$
on the summand $Y$ of $X$ corresponding to (\ref{summand})\ts.
Therefore $Y\in\A_N^{M-1}$, and we get Part (b) of Lemma \ref{L1.3}. 
\qed

\medskip\noindent
\textit{Proof of Proposition 1.4.\/}
Suppose that $|i|\com|j|\le N\ts$.
Let us prove by induction on $n=1\com2\com\ts\ldots$ that
the differences
\begin{equation}\label{diff}
F_{ij\ts|\ts N+M}^{\,(n)}-F_{ij\ts|\ts N+M-1}^{\,(n)}
\quad\text{and}\quad
C_{N+M}^{\,(n)}-C_{N+M-1}^{\,(n)}
\end{equation}
belong to the left ideal $\J_{\ts N+M}$ in $\UNM\ts$. 
Due to
(\ref{intersect}), Proposition \ref{P1.4} 
will then follow. Neither of the elements 
$F_{ij\ts|\ts N+M}^{\,(n)}$ and $C_{N+M}^{\,(n)}\,$,
nor the ideal $\J_{\ts N+M}$ depend
on the partition of the number $N+M$ into $N$ and $M$.
Hence it suffices to consider the case $M=1\ts$.
Note that according to the definition~(\ref{CN})
$$
C_{N+1}^{\,(n)}-\ts C_{N}^{\,(n)}\,=
\sum_{|k|\le N}
(\,F_{kk\ts|\ts N+1}^{\,(n)}-\ts F_{kk\ts|\ts N}^{\,(n)}\ts)
\ \,+
$$
\vspace{-8pt}
$$
F_{N+1,N+1\ts|\ts N+1}^{\,(n)}
\,+\,F_{-N-1,-N-1\ts|\ts N+1}^{\,(n)}
$$
where the last two summands belong to $\J_{\ts N+1}\ts$,
by their definition.
Therefore it suffices to consider only the first of the
differences (\ref{diff}), where $M=1\ts$. 

If $n=1$, that difference is zero.
Now suppose that $n>1$, and make the induction assumption.
Using the relation (\ref{split}),
$$
F_{ij\ts|\ts N+1}^{\,(n)}-F_{ij\ts|\ts N}^{\,(n)}
\,=\sum_{|k|\le N}\,
(-1)^{\ts\bk}\,
F_{ik}\,
(\,F_{kj\ts|\ts N+1}^{\,(n-1)}-F_{kj\ts|\ts N}^{\,(n-1)}\ts)\ \,+
$$
\vspace{-5pt}
$$
F_{i\ts,\ts N+1}\,F_{N+1,\ts j\ts|\ts N+1}^{\,(n-1)}-\ts
F_{i\ts,\ts-N-1}\,F_{-N-1,\ts j\ts|\ts N+1}^{\,(n-1)}\,\ts.
$$
At the right hand side of the last equality,
the summands corresponding to $|k|\le N$
belong to the left ideal $\J_{\ts N+1}$ by the induction assumption.
Using (\ref{FNR}) and (\ref{FN-}), 
the remainder of the right hand side is equal to the sum 
$$
(-1)^{\,\bi\,\bj}\,\,F_{N+1,\ts j\ts|\ts N+1}^{\,(n-1)}\,F_{i\ts,\ts N+1}
-\ts
(-1)^{\,\ts(\bi\,\ts+\ts1)\ts(\bj\,\ts+\ts1)}\,
F_{-N-1,\ts j\ts|\ts N+1}^{\,(n-1)}\,F_{i\ts,\ts-N-1}\ +
$$
$$
(-1)^{\,\bi\,\bj\,\,+\,1}\,
(\,1+(-1)^{\,\bi\ts\,+\,\bj\ts\,+\,n}\ts)\,\,
\de_{ij}\,\,
F_{N+1,\ts N+1\ts|\ts N+1}^{\,(n-1)}\ +
$$
$$
(-1)^{\,(\bi\,+1)\,(\bj\,+1)}\,
(\,1+(-1)^{\,\bi\ts\,+\,\bj\ts\,+\,n}\ts)\,\,
\,\de_{i,-j}\,\,
F_{-N-1,\ts N+1\ts|\ts N+1}^{\,(n-1)}\ .
$$
But in this sum, every summand evidently belongs
to the left ideal $\J_{\ts N+1}\,$.
\qed

\medskip
Let us now prove Theorem \ref{T1.5}. Firstly, we will verify
the formula (\ref{defrel}) for~the supercommutator
$[\,F_{ij}^{\ts(m)},\ts F_{kl}^{\ts(n)}\ts]$
in the algebra $\A_N\ts$. We will use

\begin{proposition*}\label{P3.1}
In\/ $\UN$ for\/
$m\com n=1\com2\ts\com\ts\ldots$ and\,
$i\com j\com k\com l=\pm1\lc\!\pm N$
$$
[\ts F_{ij\ts|\ts N}^{\,(m)}\com F_{kl\ts|\ts N}^{\,(n)}\ts]\,=\,
F_{il\ts|\ts N}^{\,(m+n-1)}\,\de_{kj}-
(-1)^{\,(\ts\bi\,+\ts\bj\ts)\ts(\ts\bk\,+\,\bl\ts)}\,
\de_{il}\,F_{kj\ts|\ts N}^{\,(m+n-1)}\,\,+
$$
$$
(-1)^{\,m\,-\ts1}\,(\,
F_{-i,l\ts|\ts N}^{\,(m+n-1)}\,\de_{-k,j}\,-
(-1)^{\,(\ts\bi\,+\ts\bj\ts)\ts(\ts\bk\,+\,\bl\ts)}\,
\de_{i,-l}\,F_{k,-j\ts|\ts N}^{\,(m+n-1)}\,)\,\,\,+
$$
$$
(-1)^{\ts\,\bj\,\bk\ts\,+\,\bj\,\bl\ts\,+\,\bk\,\bl}\,\,
\sum_{r=1}^{m-1}\,
(\,
F_{il\ts|\ts N}^{\,(n+r-1)}\,F_{kj\ts|\ts N}^{\,(m-r)}-
F_{il\ts|\ts N}^{\,(m-r)}\,F_{kj\ts|\ts N}^{\,(n+r-1)}
\,)\,\,\,+\hspace{46pt}\vspace{-6pt}
$$
$$
(-1)^{\ts\,\bj\,\bk\ts\,+\,\bj\,\bl\ts\,+\,\bk\,\bl\ts\,+\,\bk\ts\,+\,\bl}
\,\,
\sum_{r=1}^{m-1}\,(-1)^{\ts r}\,
(\,
F_{-i,l\ts|\ts N}^{\,(n+r-1)}\,F_{-k,j\ts|\ts N}^{\,(m-r)}-
F_{i,-l\ts|\ts N}^{\,(m-r)}\,F_{k,-j\ts|\ts N}^{\,(n+r-1)}
\,)\,\ts.
$$
\end{proposition*}

\begin{proof}
The formula for the supercommutator
$[\ts F_{ij\ts|\ts N}^{\,(m)}\com F_{kl\ts|\ts N}^{\,(n)}\ts]$
in the $\ZZ_2\ts$-\ts graded algebra $\UN$ displayed above
is easy to verify by using the induction on $m\ts$.
When $m=1$, this formula coincides with (\ref{FNR}).
Now make the induction assumption. Let the index $h$
run through $\pm1\lc\!\pm\!N\ts$. Then due to (\ref{split}),
$$
[\,F_{ij\ts|\ts N}^{\,(m+1)}\com F_{kl\ts|\ts N}^{\,(n)}\ts]\,=\,
\sum_h\,
(-1)^{\ts\bh}\,
[\,F_{ih}\,F_{hj\ts|\ts N}^{\,(m)}\com F_{kl\ts|\ts N}^{\,(n)}\ts]\,=
$$
$$
\sum_h\,
(-1)^{\ts\bh}\,
F_{ih}\,[\,F_{hj\ts|\ts N}^{\,(m)}\com F_{kl\ts|\ts N}^{\,(n)}\ts]\,+
\sum_h\,
(-1)^{\,\bh\ts\,+\,(\bh\ts\,+\,\bj\,)\ts(\bk\ts\,+\,\bl\,)}\,
[\,F_{ih}\com F_{kl\ts|\ts N}^{\,(n)}\ts]\,
F_{hj\ts|\ts N}^{\,(m)}
$$
$$
=\,\,
\sum_h\,
(-1)^{\ts\bh}\,
F_{ih}\,(\,
F_{hl\ts|\ts N}^{\,(m+n-1)}\,\de_{kj}-
(-1)^{\,(\ts\bh\,+\ts\bj\ts)\ts(\ts\bk\,+\,\bl\ts)}\,
\de_{hl}\,F_{kj\ts|\ts N}^{\,(m+n-1)}\,\,+
$$
$$
(-1)^{\,m\,-\ts1}\,(\,
F_{-h,l\ts|\ts N}^{\,(m+n-1)}\,\de_{-k,j}\,-
(-1)^{\,(\ts\bh\,+\ts\bj\ts)\ts(\ts\bk\,+\,\bl\ts)}\,
\de_{h,-l}\,F_{k,-j\ts|\ts N}^{\,(m+n-1)}\,)\,\,\,+
$$
$$
(-1)^{\ts\,\bj\,\bk\ts\,+\,\bj\,\bl\ts\,+\,\bk\,\bl}\,\,
\sum_{r=1}^{m-1}\,
(\,
F_{hl\ts|\ts N}^{\,(n+r-1)}\,F_{kj\ts|\ts N}^{\,(m-r)}-
F_{hl\ts|\ts N}^{\,(m-r)}\,F_{kj\ts|\ts N}^{\,(n+r-1)}
\,)\,\,\,+\hspace{46pt}\vspace{-4pt}
$$
$$
(-1)^{\ts\,\bj\,\bk\ts\,+\,\bj\,\bl\ts\,+\,\bk\,\bl\ts\,+\,\bk\ts\,+\,\bl}
\,\,
\sum_{r=1}^{m-1}\,(-1)^{\ts r}
(\,
F_{-h,l\ts|\ts N}^{\,(n+r-1)}\,F_{-k,j\ts|\ts N}^{\,(m-r)}-
F_{h,-l\ts|\ts N}^{\,(m-r)}\,F_{k,-j\ts|\ts N}^{\,(n+r-1)}
\,)\ts)\vspace{2pt}
$$
$$
+\,\,
\sum_h\,
(-1)^{\,\bh\ts\,+\,(\bh\ts\,+\,\bj\,)\ts(\bk\ts\,+\,\bl\,)}\,
(\,\ts
\de_{kh}\,F_{il\ts|\ts N}^{\,(n)}-
(-1)^{\,(\ts\bi\,+\ts\bh\ts)\ts(\ts\bk\,+\,\bl\ts)}\,
\de_{il}\,F_{kh\ts|\ts N}^{\,(n)}\,\,+\vspace{-2pt}
$$
$$
\de_{-k,h}\,F_{-i,l\ts|\ts N}^{\,(n)}-
(-1)^{\,(\ts\bi\,+\ts\bh\ts)\ts(\ts\bk\,+\,\bl\ts)}\,
\de_{i,-l}\,F_{k,-h\ts|\ts N}^{\,(n)}\ts)\,
F_{hj\ts|\ts N}^{\,(m)}\,.
$$
Here we used the
induction assumption with the index $i$ replaced by $h\ts$,
and the relation (\ref{FNR}) with the index $j$ replaced by $h\ts$.
Using the relation (\ref{split}) repeatedly,
the right hand side of the above displayed equalities equals
$$
F_{il\ts|\ts N}^{\,(m+n)}\,\de_{kj}-
(-1)^{\,(\ts\bi\,+\ts\bj\ts)\ts(\ts\bk\,+\,\bl\ts)}\,
\de_{il}\,F_{kj\ts|\ts N}^{\,(m+n)}\,\,+
$$
$$
(-1)^{\ts m}\,(\,
F_{-i,l\ts|\ts N}^{\,(m+n)}\,\de_{-k,j}\,-
(-1)^{\,(\ts\bi\,+\ts\bj\ts)\ts(\ts\bk\,+\,\bl\ts)}\,
\de_{i,-l}\,F_{k,-j\ts|\ts N}^{\,(m+n)}\,)\,\,\,+
\vspace{4pt}
$$
$$
(-1)^{\ts\,\bj\,\bk\ts\,+\,\bj\,\bl\ts\,+\,\bk\,\bl}\,\,
(\,
F_{il\ts|\ts N}^{\,(n)}\,F_{kj\ts|\ts N}^{\,(m)}\ts-\,
F_{il\ts|\ts N}^{\,(1)}\,F_{kj\ts|\ts N}^{\,(m+n-1)}
\,)\,\,\,+\vspace{4pt}
$$
$$
(-1)^{\ts\,\bj\,\bk\ts\,+\,\bj\,\bl\ts\,+\,\bk\,\bl\ts\,+\,\bk\ts\,+\,\bl}
(\ts(-1)^{\ts m+1}\,
F_{i,-l\ts|\ts N}^{\,(1)}\,F_{k,-j\ts|\ts N}^{\,(m+n-1)}-
F_{-i,l\ts|\ts N}^{\,(n)}\,F_{-k,j\ts|\ts N}^{\,(m)}
\,)\,\,\,+\vspace{-2pt}
$$
$$
(-1)^{\ts\,\bj\,\bk\ts\,+\,\bj\,\bl\ts\,+\,\bk\,\bl}\,\,
\sum_{r=1}^{m-1}\,
(\,
F_{il\ts|\ts N}^{\,(n+r)}\,F_{kj\ts|\ts N}^{\,(m-r)}-
F_{il\ts|\ts N}^{\,(m-r+1)}\,F_{kj\ts|\ts N}^{\,(n+r-1)}
\,)\,\,\,+\hspace{56pt}\vspace{-8pt}
$$
$$
(-1)^{\ts\,\bj\,\bk\ts\,+\,\bj\,\bl\ts\,+\,\bk\,\bl\ts\,+\,\bk\ts\,+\,\bl}
\,\,
\sum_{r=1}^{m-1}\,(-1)^{\ts r+1}
(\,
F_{-i,l\ts|\ts N}^{\,(n+r)}\,F_{-k,j\ts|\ts N}^{\,(m-r)}\,+\,
F_{i,-l\ts|\ts N}^{\,(m-r+1)}\,F_{k,-j\ts|\ts N}^{\,(n+r-1)}
\,)\,.
$$
But the last displayed sum can also be obtained by replacing 
$m$ by $m+1$ at right hand side of the equality in Proposition \ref{P3.1}. 
Thus we have made the induction step. 
\qed
\end{proof}

If $m>n$ then at the right hand side
of the equality in Proposition~\ref{P3.1},
the summands corresponding to the indices $r=1\lc m-n$ 
cancel in each of the two sums over $r=1\lc m-1\ts$.
In the first of the two sums, this is obvious.
To cancel these summands in the second sum, one utilises 
the relations (\ref{FN-}). Hence if $m>n$,
the summation over $r=1\lc m-1$ in Proposition \ref{P3.1}
can be replaced by the summation over $r=m-n+1\lc m-1$.
Thus if we change the running index $r$ to $m-r$,
the latter index should run through $1\lc\!\min(\ts m\com n)-1\ts$.
Using this remark, the relation (\ref{defrel}) in Theorem \ref{T1.5}
follows from Proposition \ref{P3.1}.

In the remainder of the proof of Theorem \ref{T1.5}, we will also
make use of the next proposition. For any integers $M\ge0$ and $n\ge1$ 
consider the elements
\begin{equation}\label{cf}
c_{N+M}^{\ts(n)}
\quad\textrm{and}\quad
f_{ij\ts|\ts N+M}^{\,(n)}
\end{equation}
of the algebra $\S(\qNM)^{\,\qM}$,
see Proposition \ref{P2.5}. 
Fix any positive integer $s\ts$.

\begin{proposition*}\label{P3.2}
Take the elements 
$f_{ij\ts|\ts N+M}^{\,(n)}$
where 
\begin{equation}\label{tricon}
1\le n\le s\ts,\ \ 
1\le i\le N\ts,\ \ 
1\le |j|\le N\ts.
\end{equation}
Along with these elements, take
the elements $c_{N+M}^{\ts(n)}$ where $1\le n\le s$ and
$n$ is odd. For any sufficiently large number $M$,
all these elements
are algebraically independent
in the supercommutative algebra $\S(\qNM)\ts$.
\end{proposition*}

\begin{proof}
We will use arguments from \cite[Subsection 2.11]{MO}.
By the Poincar\'e\,-Birkhoff\,-Witt theorem 
for Lie superalgebras, the elements 
$$
F_{kl}=f_{kl\ts|\ts N+M}^{\,(1)}
$$ where
$k=1\lc N+M$ and $l=\pm1\lc\!\pm\!N+M\ts$, are 
free generators of the supercommutative algebra $\S(\qNM)\ts$.
Let $\X_{\ts s}$ be the quotient algebra of $\S(\qNM)\ts$,
defined by imposing the following relations on these free generators.

For every triple $(\ts i\com j\com n)$ satisfying the 
conditions (\ref{tricon}), choose a subset
$$
\O_{ij}^{\ts(n)}\subset\{\ts N+1\com\ts N+2\com\,\ldots\,\ts\}
$$
of cardinality $n-1$ in such a way, that all these subsets are disjoint.
Let $M$ be so large that all these subsets are
contained in $\{\ts N+1\lc N+M-s\ts\}$. If
$$
\O_{ij}^{\ts(n)}=\{\ts l_1\lc l_{n-1}\ts\}\,,
$$
then put
$$
F_{\ts i\ts l_1}=F_{\,l_1l_2}=\ldots=F_{\,l_{n-2}\ts l_{n-1}}=1\ts.
$$
Denote by $x_{ij}^{\ts(n)}$
the image of the element $F_{\,l_{n-1}\ts j}\in\S(\qNM)$ in
the algebra $\X_{\ts s}\ts$.
Having done this for every triple $(\ts i\com j\com n)$ satisfying
the conditions 
(\ref{tricon}), for every $r=1\lc s$ denote by $x_r$ the image in
$\X_{\ts s}$ of the element 
$$
F_{\,N+M-s+r\ts,\ts N+M-s+r}\in\S(\qNM)\ts.
$$
Finally, put $F_{kl}=0$ if
$k>0$ and $(k\com l)$ is not one of the pairs
$$
(\ts i\com l_1)
\ts\com\ts
(\ts l_1\com l_2)
\,\lc
(\ts l_{n-2}\com l_{n-1})
\ts\com\ts
(\ts l_{n-1}\com j)
$$
for any triple $(\ts i\com j\com n)$ satisfying
(\ref{tricon}), and not one of the pairs
$$
(N\!+\!M\!-\!s\!+\!r\com N\!+\!M\!-\!s\!+\!r)
\ \ \textrm{where}\ \ 
r=1\lc s\,.
$$ 

The elements $x_1\lc x_s\,$ and the elements
$x_{ij}^{\ts(n)}$ for all the triples $(i\com j\com n)$
satisfying (\ref{tricon}), are free generators of the
algebra $\X_{\ts s}\ts$. For any of these triples,
the image in $\X_{\ts s}$ of $$f_{ij\ts|\ts N+M}^{\,(n)}\in\S(\qNM)$$
equals $x_{ij}^{\ts(n)}$ plus a certain linear combination
of products of the elements $x_{kl}^{\ts(m)}$
where $1\le m<n$.
For any odd $n\ts$, the image in $\X_{\ts s}$
of the element $$c_{N+M}^{\,(n)}\in\S(\qNM)$$ 
equals
$$2^{\ts n}(\ts x_1^n+\ldots+x_s^n\ts)$$ plus a linear combination
of products of elements $x_{kl}^{\ts(m)}$ where $1\le m\le n\ts$. 
Hence all these images are algebraically independent in 
the quotient $\X_{\ts s}$ of the supercommutative algebra $\S(\qNM)\ts$.
\qed
\end{proof}

Let us show that the associative
algebra $\A_N$ is generated by the elements
$C^{\ts(1)},C^{\ts(3)},\,\ldots\,$ and
$F_{ij}^{\ts(1)},F_{ij}^{\ts(2)},\,\ldots\,\,$.
Take any element $Z\in\A_N\ts$,
and consider its canonical 
image $Z_M=\pi_M(Z)\in\AMN$ for any $M\ge0\ts$.
By Theorem \ref{T1.1}, the element $Z_M$ is a linear combination of
the products of the elements $C_{N+M}^{\,(n)}$ 
where $n=1\com3\com\,\ldots$ and of the elements
$F_{ij\ts|\ts N+M}^{\,(n)}$ where $n=1\com2\com\,\ldots\,\,$ whereas
$i=1\lc N$ and $j=\pm1\lc\!\pm\!N\ts$;
see (\ref{FN-}).
Choose any linear ordering 
of all these elements.
Applying Proposition \ref{P3.1} 
to the algebra $\UNM$ instead of $\UN\ts$, we will assume
that any of the products in the linear combination $Z_M$ 
is an ordered monomial in these elements.
If $\bi+\bj=1$, then the element
$F_{ij\ts|\ts N+M}^{\,(n)}$
may appear in any of these monomials only with 
the degree $1$.

We will assume that for any $M\ge1$, the map
$\al_M$ preserves the ordering. Then for every monomial $Y_M$
appearing in the linear combination $Z_M\ts$, the monomial
$\al_M(\ts Y_M)$ may appear in the linear combination $\al_M(Z_M)=Z_{M-1}$.

The filtration degrees of all the elements 
$Z_0\com Z_1\com Z_2\com\,\ldots$
are bounded from above. Hence for any factor $C_{N+M}^{\,(n)}$
or $F_{ij\ts|\ts N+M}^{\,(n)}$ of the monomials appearing in
linear combination $Z_M\ts$, we have $n\le s$ 
for a certain integer $s$ which does not depend on $M$. 
Then for a sufficiently large number $M$,
the coefficients of the monomials appearing in the 
linear combinations $Z_M\com Z_{M+1}\com Z_{M+2}\com\,\ldots$
are determined uniquely. The uniqueness follows from Proposition \ref{P3.2}.

Now fix a sufficiently large number $M$, as above.
Let $Y_M$ be any
monomial appearing in the linear combination $Z_M\ts$,
say with a coefficient $z\in\CC\ts$. We assume that $z\neq0\ts$.
Then for any integer $L>M$, the linear combination $Z_L$
contains the summand
$z\ts Y_{L}$ where $Y_L$ is a monomial and 
$$
(\ts\al_{M+1}\circ\ldots\circ\al_L\ts)(\ts Y_L)=Y_M\ts.
$$
For any non-negative integer $L\le M$, define
$
Y_L=(\ts\al_L\circ\ldots\circ\al_M\ts)(\ts Y_M)\ts.
$
The sequence $Y_0\com Y_1\com Y_2\com\,\ldots$
determines an element $Y\in\A_N$,
which is a monomial in
$C^{\ts(1)},C^{\ts(3)},\,\ldots\,\,$ and
$F_{ij}^{\ts(1)},F_{ij}^{\ts(2)},\,\ldots\,\,$.
Here $1\le i\le N$ and $1\le|j|\le N$.
The element $Z\in\A_N$ is then a sum of
the products of the form $z\ts Y$.
This sum is finite, because
any such product corresponds to a summand
$z\ts Y_M$ in the linear combination $Z_M\ts$.
Thus we have proved Part (a) of Theorem \ref{T1.5}.

Let us now prove Parts (b) and (c).
By definition, the algebra $\A_N$ comes with an ascending
$\ZZ\ts$-\ts filtration, such that the generators 
$C^{\ts(n)}$ and $F_{ij}^{\ts(n)}$ of
$\A_N$ have the degree $n\ts$. 
Denote the corresponding $\ZZ\ts$-\ts graded algebra
by $\,{\rm{gr}}\,\A_N\ts$. Let $c^{\ts(n)}$ and $f_{ij}^{\ts(n)}$
be the generators of the algebra $\,{\rm{gr}}\,\A_N$
corresponding to $\,C^{\ts(n)}$ and $F_{ij}^{\ts(n)}$.
We always assume that the index $n$ in $C^{\ts(n)}\,$ and $c^{\ts(n)}$
is odd. We also assume that $|i|\com|j|\le N$ in
$F_{ij}^{\ts(n)}$ and $f_{ij}^{\ts(n)}$. It follows from (\ref{F-}) that
for any $n=1\com2\com\,\ldots$
$$
f_{-i,-j}^{\,(n)}\ts=
(-1)^{\,n-1}\,
f_{ij}^{\,(n)}.
$$

The algebra ${\rm{gr}}\,\A_N\ts$
also inherits from $\A_N$ a $\ZZ_2\ts$-\ts gradation, such that
$$
\deg\ts c^{\,(n)}=0
\quad\textrm{and}\quad
\deg\ts f_{ij}^{\,(n)}=\bi\ts+\bj\,,
$$
see (\ref{grad}).
The relation (\ref{defrel}) demonstrates that the $\ZZ_2\ts$-\ts graded 
algebra $\,{\rm{gr}}\,\A_N\ts$ is supercommutative. To complete the proof
of Theorem \ref{T1.5}, it suffices to show that the elements 
$c^{\ts(1)},c^{\ts(3)},\,\ldots\,$ together with the elements
$f_{ij}^{\ts(1)},f_{ij}^{\ts(2)},\,\ldots$ where
$i=1\lc N$ and $j=\pm1\lc\!\pm\!N\ts$, are
algebraically independent in the
supercommutative algebra $\,{\rm{gr}}\,\A_N\ts$.

The algebra $\,{\rm{gr}}\,\A_N\ts$
can also be obtained as an inverse limit 
of the sequence
of the supercommutative algebras
$\S(\qNM)^{\,\qM}$ where $M=0\com1\com2\com\,\ldots\ts\,$.
The limit is taken in the category of $\ZZ\ts$-\ts graded algebras.
We assume that if $M=0$, then $\S(\qNM)^{\,\qM}=\S(\qN)\ts$.
The definition of the surjective homomorphism
$$
\S(\qNM)^{\,\qM}\to\,\S(\q_{\ts N+M-1})^{\,\q_{M-1}}
$$
for any $M\ge1$ is similar to the definition of the surjective
homomorphism $\al_M:\AMN\to\A_N^{M-1}$, see Lemma \ref{L1.3}.
Here we omit the details, but notice that the elements
$c^{\ts(n)}$ and $f_{ij}^{\ts(n)}$ of $\,{\rm{gr}}\,\A_N\ts$
correspond to the sequences of elements (\ref{cf})
of the algebras $\S(\qNM)^{\,\qM}$
where $M=0\com1\com2\com\,\ldots\,\,$.
Proposition \ref{P3.2} now guarantees
the algebraic independence of
the elements 
$c^{\ts(1)},c^{\ts(3)},\,\ldots\,$ together with the elements
$f_{ij}^{\ts(1)},f_{ij}^{\ts(2)},\,\ldots$ where
$i=1\lc N$ and $j=\pm1\lc\!\pm\!N\ts$.

\medskip\noindent
\textit{Proof of Proposition 1.6.\/}
Under the correspondence
$F_{ij}^{\ts(n)}\mapsto(-1)^{\,\bi}\ts\,T_{ji}^{\ts(n)}$,
the collection of relations (\ref{F-}) in the algebra $\B_N$
for the indices $n=1\com2\com\,\ldots$
corresponds to the equality (\ref{Fu-}) in $\YN\ts[[\ts u^{\ts-1}]]\ts$.
Put $T_{ij}^{\ts(0)}=\de_{ij}$ for any $i\com j=1\lc\!\pm\!N$.
Using ({\ref{ai}), the relation (\ref{defrel}) in $\B_N$
then corresponds to
$$
(-1)^{\ts\,\bi\,\bj\ts\,+\,\bi\,\bk\ts\,+\,\bj\,\bk}\,\,
[\,\ts T_{ji}^{\ts(m)},\ts T_{lk}^{\ts(n)}\ts]
\,=\,
\vspace{-4pt}
$$
$$
\sum_{r=0}^{m-1}\,
(\,
T_{jk}^{\ts(m+n-r-1)}\,T_{li}^{\ts(r)}-\,
T_{jk}^{\ts(r)}\,T_{li}^{\ts(m+n-r-1)}
\,)\,\,+\,\,
(-1)^{\ts\,\bi\ts\,+\,\bj\,\ts+\ts1}
\,\,\times\hspace{20pt}\vspace{-6pt}
$$
\begin{equation}\label{jikl}
\hspace{50pt}
\sum_{r=0}^{m-1}
(-1)^{\,m\,+\,r}\,
(\,
T_{-j,k}^{\ts(m+n-r-1)}\,T_{-l,i}^{\ts(r)}-
T_{j,-k}^{\ts(r)}\,T_{l,-i}^{\ts(m+n-r-1)}
\,)\ .
\vspace{2pt}
\end{equation}
Here we also used a remark 
on the summation over $r=1\lc m-1\ts$ similar to that made immediately after
the proof of Proposition \ref{P3.1}.

Put $T_{ij}^{\ts(-1)}=0$ for any $i\com j=1\lc\!\pm N$.
The collection of relations (\ref{jikl}) for
$m\com n=1\com2\com\,\ldots$
is equivalent to the collection of relations
$$
(-1)^{\ts\,\bi\,\bj\ts\,+\,\bi\,\bk\ts\,+\,\bj\,\bk}\,\,
(\,[\,\ts T_{ji}^{\ts(m+1)},\ts T_{lk}^{\ts(n-1)}\ts]-
[\,\ts T_{ji}^{\ts(m-1)},\ts T_{lk}^{\ts(n+1)}\ts]\,)
\,=\,
$$
$$
T_{jk}^{\ts(n-1)}\,T_{li}^{\ts(m)}-\,
T_{jk}^{\ts(m)}\,T_{li}^{\ts(n-1)}+\,
T_{jk}^{\ts(n)}\,T_{li}^{\ts(m-1)}-\,
T_{jk}^{\ts(m-1)}\,T_{li}^{\ts(n)}\,+\,
(-1)^{\ts\,\bi\ts\,+\,\bj}
\,\,\,\times\vspace{2pt}
$$
\begin{equation}\label{jikl-}
(\,
T_{-j,k}^{\ts(n-1)}\,T_{-l,i}^{\ts(m)}-
T_{j,-k}^{\ts(m)}\,T_{l,-i}^{\ts(n-1)}-
T_{-j,k}^{\ts(n)}\,T_{-l,i}^{\ts(m-1)}+
T_{j,-k}^{\ts(m-1)}\,T_{l,-i}^{\ts(n)}
\,)
\vspace{4pt}
\end{equation}
for $m\com n=0\com1\com2\com\,\ldots\ $.
Multiplying the relation (\ref{jikl-}) by $x^{\ts1-m}\,y^{\ts1-n}$
and taking the sum of resulting relations
over $m\com n=0\com1\com2\com\,\ldots$ we 
get the relation
$$
\vspace{2pt}
(\ts x^2-y^2\ts)\cdot
[\,T_{ji}(x)\ts,T_{lk}(y)\ts]
\cdot
(-1)^{\ts\,\bi\,\bj\ts\,+\,\bi\,\bk\ts\,+\,\bj\,\bk}
\,=\,
$$
$$
(x+y)
\cdot
(\ts T_{jk}(y)\,T_{li}(x)-T_{jk}(x)\,T_{li}(y)\ts)\ +
\vspace{4pt}
$$
\begin{equation}\label{wrongrel}
(x-y)
\cdot
(\ts T_{-j,\ts k}(y)\,T_{-l,\ts i}(x)-T_{j,\ts -k}(x)\,T_{l,\ts -i}(y) \ts)
\cdot
(-1)^{\ts\bi\ts\,+\ts\bj}\,.
\vspace{4pt}
\end{equation}
Using (\ref{scom}), we can rewrite the left hand
side of the relation (\ref{wrongrel}) as
$$
(\ts y^2-x^2\ts)\cdot
[\,T_{lk}(y)\ts,T_{ji}(x)\ts]
\cdot
(-1)^{\ts\,\bi\,\bj\ts\,+\,\bi\,\bl\ts\,+\,\bj\,\bl}\,.
$$
Replacing in the resulting relation the indices $i\com j\com k\com l$
and the parameters $x\com y$ by $l\com k\com j\com i$ and $y\com x$
respectively, we obtain exactly the relation (\ref{yangrel})\ts. Thus
the defining relations of the subalgebra $\B_N\subset\A_N$ 
correspond to the defining relations of the algebra $\YN\ts$,
see Theorem \ref{T1.5}.
\qed



\end{document}